\documentclass[12pt,reqno]{article}

\usepackage{amssymb}
\usepackage{enumerate}
\usepackage{bm}
\usepackage{lscape}
\usepackage{booktabs,caption}
\usepackage{multirow}
\usepackage{graphicx}
\usepackage{amsmath,amsthm}
\usepackage{a4wide}
\usepackage{version}
\usepackage{parskip}
\usepackage{euscript}
\usepackage{mathrsfs}
\usepackage{tikz}
\usetikzlibrary{patterns}
\usetikzlibrary{arrows,decorations.markings, matrix}
\usepackage{mathtools}
\usepackage{hyperref}
\usepackage{comment}

\usepackage[all]{xy}

\numberwithin{equation}{section}

\theoremstyle{definition}

\newcommand{\Gm}{\mathbb{G}_m}
\newcommand{\avoids}{\mathcal{B}}
\newcommand{\Yon}{Y}

\newcommand{\Z}{\mathbb{Z}}

\newcommand{\WW}{{\mathcal W}}

\newcommand{\MM}{{\mathcal M}}

\renewcommand{\AA}{{\mathcal A}}

\newcommand{\VV}{{\mathcal V}}

\newcommand{\LL}{{\mathcal L}}
\renewcommand{\L}{{\mathbb{L}}}
\newcommand{\YY}{{\mathcal Y}}

\newcommand{\Hom}{\operatorname{Hom}}

\newcommand{\End}{\operatorname{End}}

\newcommand{\C}{{\mathbb C}}
\newcommand{\R}{{\mathbb R}}

\renewcommand{\mod}{\operatorname{mod}}

\newcommand{\OO}{{\mathcal O}}

\newcommand{\Coh}{\operatorname{coh}}

\newcommand{\Spec}{\operatorname{Spec}}

\renewcommand{\P}{{\mathbb P}}

\newcommand{\A}{{\mathbb A}}

\newcommand{\perf}{\operatorname{perf}}

\newcommand{\m}{\mathfrak{m}}

\newcommand{\OP}{\operatorname}

\makeatletter
\def\thm@space@setup{%
  \thm@preskip=\parskip \thm@postskip=0pt
}
\makeatother

\tikzset{->-/.style={decoration={
              markings,
              mark=at position .5 with {\arrow{>}}},postaction={decorate}}}

\tikzset{-<-/.style={decoration={
              markings,
              mark=at position .5 with {\arrow{<}}},postaction={decorate}}}

\newcommand{\sslash}{\mathbin{/\mkern-6mu/}}
\usetikzlibrary{arrows.meta}
\title{Noncommutative crepant resolutions of \(cA_n\) \\ singularities via Fukaya categories}
\author{Jonny Evans\hspace{2cm} Yank\i\ Lekili}
\renewcommand{\theparagraph}{\S\arabic{section}.\arabic{paragraph}\noindent}
\newcommand{\pg}{\paragraph{\hspace{-0.6cm}}}
\newcommand{\Addresses}{{
  \bigskip
  \footnotesize
  J.~D.~Evans, \textsc{University of Lancaster}\par\nopagebreak
  \texttt{j.d.evans@lancaster.ac.uk}

  \medskip

  Y.~Lekili, \textsc{Imperial College London}\par\nopagebreak
  \texttt{y.lekili@imperial.ac.uk}

}}

\date{} 

\begin{document}

\maketitle 

\begin{abstract} We compute the wrapped Fukaya category
  $\mathcal{W}(T^*S^1, D)$ of a cylinder relative to a divisor
  $D= \{p_0,\ldots, p_n\}$ of $n+1$ points, proving a mirror
  equivalence with the category of perfect complexes on a
  crepant resolution (over $k[\![t_0,\ldots, t_n]\!]$) of the
  singularity $uv=t_0t_1\ldots t_n$. Upon making the base-change
  $t_i= f_i(x,y)$, we obtain the derived category of any crepant
  resolution of the $cA_{n}$ singularity given by the equation
  $uv= f_0\ldots f_n$. These categories inherit braid group
  actions via the action on $\mathcal{W}(T^*S^1,D)$ of the
  mapping class group of $T^*S^1$ fixing $D$. We also give geometric models for the derived contraction algebras associated to
  a \(cA_n\) singularity in terms of the relative Fukaya category of the
  disc.
\end{abstract}

\section{Introduction}

\pg Consider the Fukaya category of a point with coefficients in
a ring \(R\). Before taking the triangulated envelope, there is
only one object: the point itself, with endomorphism algebra
\(R\). If \(R\) is not a field then there are non-invertible
non-zero endomorphisms which allow us to construct new twisted
complexes in the derived Fukaya category. Via the Yoneda
embedding, we can think of the derived Fukaya category of a
point with coefficients in \(R\) as \(\perf(R)\). We can think
of this as the world's lousiest \(A\)-model mirror to
\(\Spec R\). It is lousy in the precise sense that symplectic
geometry has given us absolutely no information here: all of the
interesting information is contained in the coefficient
ring. The moral of the current paper is that there is a whole
spectrum of ways we can get at a single triangulated
\(A_\infty\)-category by combining symplectic manifolds with
coefficient rings. We work out in detail some examples where the
symplectic manifold is a 2-dimensional cylinder.

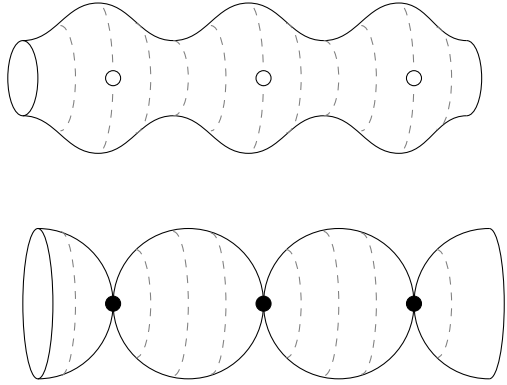
\begin{figure}[htb]
  \begin{center}
    \begin{tikzpicture}
      \draw (-3.2,0.5) to[out=0,in=180] (-2.2,1) to[out=0,in=180] (-1.2,0.5) to[out=0,in=180] (-0.2,1) to[out=0,in=180] (0.8,0.5) to[out=0,in=180] (1.8,1) to[out=0,in=180] (2.7,0.5);
      \draw (-3.2,-0.5) to[out=0,in=180] (-2.2,-1) to[out=0,in=180] (-1.2,-0.5) to[out=0,in=180] (-0.2,-1) to[out=0,in=180] (0.8,-0.5) to[out=0,in=180] (1.8,-1) to[out=0,in=180] (2.7,-0.5);
      \draw (-3.2,0) circle [x radius = 0.2cm, y radius = 0.5cm];
      \draw (2.7,0.5) arc [x radius = 0.2cm, y radius = 0.5cm, start angle=90, end angle = -90];
      \draw[dashed,gray] (-2.7,0.7) arc [x radius=0.2, y radius=0.7, start angle=90, end angle = -90];
      \draw[dashed,gray] (-1.7,0.73) arc [x radius=0.2, y radius=0.73, start angle=90, end angle = -90];
      \draw[dashed,gray] (-1.2,0.5) arc [x radius=0.2, y radius=0.5, start angle=90, end angle = -90];
      \draw[dashed,gray] (-0.7,0.7) arc [x radius=0.2, y radius=0.7, start angle=90, end angle = -90];
      \draw[dashed,gray] (0.3,0.73) arc [x radius=0.2, y radius=0.73, start angle=90, end angle = -90];
      \draw[dashed,gray] (0.8,0.5) arc [x radius=0.2, y radius=0.5, start angle=90, end angle = -90];
      \draw[dashed,gray] (1.3,0.7) arc [x radius=0.2, y radius=0.7, start angle=90, end angle = -90];
      \draw[dashed,gray] (2.3,0.68) arc [x radius=0.2, y radius=0.68, start angle=90, end angle = -90];
      \draw[dashed,gray] (-2,0) arc [x radius=0.2, y radius=1, start angle=0, end angle=90];
      \draw[dashed,gray] (-2,0) arc [x radius=0.2, y radius=1, start angle=0, end angle=-90];
      \draw[dashed,gray] (0,0) arc [x radius=0.2, y radius=1, start angle=0, end angle=90];
      \draw[dashed,gray] (0,0) arc [x radius=0.2, y radius=1, start angle=0, end angle=-90];
      \draw[dashed,gray] (2,0) arc [x radius=0.2, y radius=1, start angle=0, end angle=90];
      \draw[dashed,gray] (2,0) arc [x radius=0.2, y radius=1, start angle=0, end angle=-90];
      \filldraw[fill=white,draw=black] (-2,0) circle [radius=0.1cm];
      \filldraw[fill=white,draw=black] (0,0) circle [radius=0.1cm];
      \filldraw[fill=white,draw=black] (2,0) circle [radius=0.1cm];
      \begin{scope}[shift={(0,-3)}]
        \draw (-3,0) circle [x radius = 0.2cm, y radius = 1cm];
        \draw (3,1) arc [x radius = 0.2cm, y radius = 1cm, start angle=90, end angle = -90];
        \draw (-3,1) arc [radius=1cm, start angle = 90, end angle = -90];
        \draw (-1,0) circle [radius=1];
        \draw (1,0) circle [radius=1];
        \draw (3,1) arc [radius=1, start angle=90, end angle=270];        
        \draw[dashed,gray] (-2.7,0.96) arc [x radius=0.2, y radius=0.96, start angle=90, end angle = -90];
        \draw[dashed,gray] (-1.7,0.72) arc [x radius=0.2, y radius=0.72, start angle=90, end angle = -90];
        \draw[dashed,gray] (-1.2,0.96) arc [x radius=0.2, y radius=0.96, start angle=90, end angle = -90];
        \draw[dashed,gray] (-0.7,0.94) arc [x radius=0.2, y radius=0.94, start angle=90, end angle = -90];
        \draw[dashed,gray] (0.3,0.72) arc [x radius=0.2, y radius=0.72, start angle=90, end angle = -90];
        \draw[dashed,gray] (0.8,0.96) arc [x radius=0.2, y radius=0.96, start angle=90, end angle = -90];
        \draw[dashed,gray] (1.3,0.94) arc [x radius=0.2, y radius=0.94, start angle=90, end angle = -90];
        \draw[dashed,gray] (2.3,0.72) arc [x radius=0.2, y radius=0.72, start angle=90, end angle = -90];
        \filldraw[fill=black,draw=black] (-2,0) circle [radius=0.1cm];
        \filldraw[fill=black,draw=black] (0,0) circle [radius=0.1cm];
        \filldraw[fill=black,draw=black] (2,0) circle [radius=0.1cm];
      \end{scope}
    \end{tikzpicture}
    \caption{A punctured cylinder \(T^*S^1\setminus D\) and a
      nodal curve \(C_{n+1}\). Both are equipped with dual
      Lagrangian torus fibrations---the fibres are the dashed
      curves. The fibres above are dual to those below in the
      sense of having reciprocal radii; the noncompact fibres
      (``infinite radius'') through the punctures are dual to
      the nodes (``zero radius'').}
     \label{fig:syz}
  \end{center}
\end{figure}

\pg \label{para:lekpol} The starting point for these examples is
the mirror symmetry result proved in \cite{LP} between (on the
A-side) \(T^*S^1\) with a collection \(D\) of \(n+1\) punctures
and (on the B-side) a certain reducible curve \(C_{n+1}\) with
\(n+1\) nodes. The two sides of the mirror, together with dual
Lagrangian torus fibrations are shown in Figure \ref{fig:syz}
(the noncompact fibres on the A-side are dual to the point-like
fibres on the B-side). The precise statement of mirror symmetry
identifies the wrapped Fukaya category of Lagrangian branes
avoiding the punctures with the derived category of perfect
complexes on the nodal curve.

\pg Consider the versal deformation \(\{uv=t_0\cdots t_n\}\) of
an \(A_{n}\)-curve singularity; this admits a crepant resolution
\(\YY\) with a morphism to \(\Spec k[t_0,\ldots,t_n]\) whose
central fibre is \(C_{n+1}\). The B-model in our main example
will be \(\YY\). To build an A-model mirror to this, we need to
find a Fukaya category which is linear over
\(R=k[t_0,\ldots,t_n]\) and which specialises to the Fukaya
category of the \((n+1)\)-punctured cylinder when the
\(t\)-variables are set equal to zero. We therefore use \(R\) as
the coefficient ring\footnote{to get \(R\)-linearity.}  for
Floer theory on \(T^*S^1\) and work relative to \(D\), using
intersections with \(D\) to weight polygons contributing to the
Floer \(A_\infty\)-operations.\footnote{to get the deformation.}
We will further base-change coefficient rings to find mirrors to
non-versal deformations.

\pg Here is the general setting. Let \(\Sigma\) be a two-dimensional Liouville manifold (non-compact surface), equipped with a choice of grading data
(line field), and let
\(D=\{z_0,\ldots,z_n\}\subset \Sigma\) be a finite set of marked
points. Fix a field \(k\), let \(n=|D|-1\), and let
\(R:=k[t_0,\ldots,t_n]\). We consider the following wrapped
Fukaya category of \(\Sigma\) relative to \(D\):
\begin{itemize}
\item The objects are properly-immersed, exact, graded
  Lagrangian branes in \(\Sigma\) avoiding the marked points
  \(D\) and asymptotic to conical Lagrangians near the ends of
  \(\Sigma\). The brane-data comprises a choice of orientation,
  relative spin-structure, grading, and local system.
\item The hom-spaces are given by wrapped intersections (see
  \cite{AbouzaidSeidel} or {\cite[Appendix B]{EkLe}}).
\item The \(A_\infty\)-operations are given by counting
  holomorphic polygons with boundaries on (wrapped) Lagrangians,
  but each polygon \(P\) contributes to the corresponding
  operation with a weight of
  \(\prod_{i=0}^nt_i^{\OP{mult}(P,z_i)}\in R\).
\item Finally, we take the split-closed triangulated envelope to
  get an \(R\)-linear triangulated \(A_\infty\)-category which
  we will write as \(\WW(\Sigma,D)\).
\end{itemize}

\pg We will frequently change our coefficient ring \(R\). If
\(S\) is an \(R\)-algebra (i.e. a ring with a morphism
\(R\to S\)) then we will write \(\WW(\Sigma,D)\otimes_R S\) for
the corresponding \(S\)-linear \(A_\infty\)-category where all
hom-spaces are tensored with \(S\).

\pg Relative Fukaya categories have played an important role in
Floer theory starting with Seidel's paper on mirror symmetry for
the quartic surface \cite{SeidelQuartic}, and the idea of
deforming Floer cohomology by weighting operations according to
how many times a polygon passes through a point goes back to
Ozsv\'{a}th and Szabo \cite{OzSz} in their work on Heegaard
Floer homology. For a detailed exposition of Fukaya categories
in the exact setting, see \cite{SeidelBook}; for wrapped
categories in general, see \cite{AbouzaidSeidel} or
\cite[Appendix B]{EkLe}, but for a very explicit model of the
wrapped Fukaya category of a surface, see {\cite{Bocklandt} and
  \cite[Section 3.3]{HKK}}. For relative (wrapped) Fukaya categories see
\cite{mamaev, PerutzSheridan,SheridanVersality} and for a very similar
example of a relative Fukaya category of a surface, see
\cite{LekiliPolishchukM1n}, and for a version with an arithmetic flavour see \cite{LekiliTreumann}.  

\paragraph{Main Theorem.}\label{thm:main}
We will focus on the specific case where \(\Sigma\) is the
cotangent bundle \(T^*S^1\) with its canonical exact symplectic
form and the line field given by cotangent fibres. We will pick
a collection of Lagrangian arcs \(L_0,\ldots,L_n\) as shown in
Figure \ref{fig:sigma}. Let \(S\) be an \(R\)-algebra. We will
prove the following results:

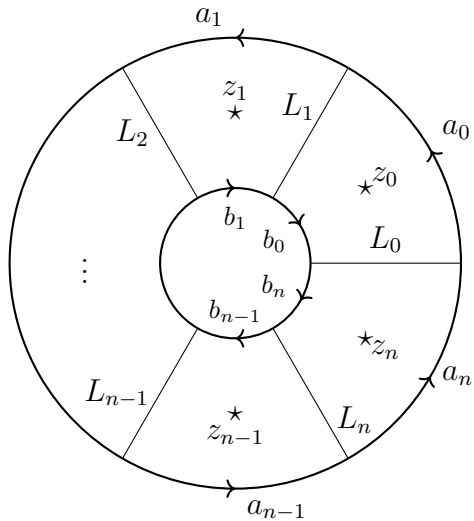
\begin{figure}[htb]
  \begin{center}
    \begin{tikzpicture}
      \draw[thick,->-] (00:3) arc[radius = 3,start angle=0,end angle=60] node [midway,above right] {\(a_{0}\)};
      \draw[thick,->-] (60:3) arc[radius = 3,start angle=60,end angle=120] node [midway,above left] {\(a_{1}\)};
      \draw[thick] (120:3) arc[radius = 3,start angle=120,end angle=240];
      \draw[thick,->-] (240:3) arc[radius = 3,start angle=240,end angle=300] node [midway,below right] {\(a_{n-1}\)};
      \draw[thick,->-] (300:3) arc[radius = 3,start angle=300,end angle=360] node [midway,right] {\(a_{n}\)};
      \draw[thick,->-] (60:1) arc[radius = 1,start angle=60, end angle=0];
      \draw[thick,->-] (0:1) arc[radius = 1,start angle=0, end angle=-60];
      \draw[thick,->-] (-60:1) arc[radius = 1,start angle=-60, end angle=-120];
      \draw[thick] (-120:1) arc[radius = 1,start angle=-120, end angle=-240];
      \draw[thick,->-] (120:1) arc[radius = 1,start angle=120, end angle=60];
      \node at (30:0.6) {\footnotesize\(b_{0}\)};
      \node at (-30:0.6) {\footnotesize\(b_{n}\)};
      \node at (-90:0.65) {\footnotesize\(b_{n-1}\)};
      \node at (90:0.6) {\footnotesize\(b_{1}\)};
      \node at (30:2cm) {\(\star\)};
      \node at (90:2cm) {\(\star\)};
      \node at (270:2cm) {\(\star\)};
      \node at (330:2cm) {\(\star\)};
      \node at (60:2.4cm) [left] {\(L_1\)};
      \node at (120:2cm) [left] {\(L_2\)};
      \node at (240:2cm) [left] {\(L_{n-1}\)};
      \node at (300:2.4cm) [right] {\(L_{n}\)};
      \node at (0:2cm) [above] {\(L_0\)};
      \node at (30:2.3cm) {\(z_0\)};
      \node at (90:2.3cm) {\(z_1\)};
      \node at (270:2.3cm) {\(z_{n-1}\)};
      \node at (330:2.3cm) {\(z_n\)};
      \node at (180:2cm) {\(\vdots\)};
      \draw (0:1cm) -- (0:3cm);
      \draw (60:1cm) -- (60:3cm);
      \draw (120:1cm) -- (120:3cm);
      \draw (240:1cm) -- (240:3cm);
      \draw (300:1cm) -- (300:3cm);
    \end{tikzpicture}
  \end{center}
  \caption{The surface \(T^*S^1\) together with its Lagrangian
    arcs \(L_0,\ldots,L_n\), marked points \(z_0,\ldots,z_n\)
    and some of the Reeb chords \(a_{i}\) and \(b_{i}\).}
  \label{fig:sigma}
\end{figure}

\begin{enumerate}[A.]
\item {\it The endomorphism \(A_\infty\)-algebra of
    \(\bigoplus_{i=0}^nL_i\) in \(\WW(T^*S^1,D)\otimes_R S\) is
    quasi-isomorphic to the algebra \(\AA(T^*S^1,D)\otimes_R S\)
    where \(\AA(T^*S^1,D)\) is defined in \ref{dfn:alg}
    below. This algebra is supported in degree zero, and hence
    has no nontrivial higher products.}  (See Section
  \ref{sct:hf_comp}.)
\item {\it Let \(\LL \subset\WW(T^*S^1,D)\) denote the
    subcategory split-generated by the Lagrangian arcs
    \(L_0,\ldots,L_n\). Then \(\LL\, \otimes_R S\) is preserved
    by the action of the mapping class group
    \(\Gamma(T^*S^1,D)\) of compactly-supported graded
    symplectomorphisms of \(T^*S^1\) fixing \(D\) pointwise.}
  (See Section \ref{sct:group_action}.)
\end{enumerate}

\paragraph{Remarks.} (i) In Appendix \ref{sct:generation}, we
will show that the arcs split-generate the category
\(\WW(T^*S^1,D)\otimes_R \bar{R}\) where \(\bar{R}\) is the
completion \(k[\![t_0,\ldots,t_n]\!]\). We expect that the arcs generate \(\WW(T^*S^1,D)\) itself, and this is confirmed in the forthcoming work of Mamaev.

(ii) We will prove something slightly more general than
\ref{thm:main}(B) which gives quasi-equivalences for
symplectomorphisms which permute the points of \(D\). For some
choices of \(R\)-algebra \(S\), these will be autoequivalences
of \(\LL\). See \ref{pg:group_action} for details.

(iii) By construction the algebra $\mathcal{A}(T^*S^1,D)$ is
linear over $R$ but, in fact, it turns out that it has a bigger
center given by $R[u,v]/(uv - t_0t_1\ldots t_n)$. We expect that
the autoequivalences given in \ref{thm:main}(B) are linear over
this bigger ring (not just linear over $R$). The main reason to
expect this is that the additional variables $u$ and $v$ come
from Hochschild cohomology classes of $\mathcal{A}(T^*S^1,D)$
associated with the infinite ends of $T^*S^1$, whereas our
autoequivalences are induced by compactly supported
symplectomorphisms.

\paragraph{Mirror symmetry
  interpretation.}\label{pg:mirror_interp}
Theorem \ref{thm:main}(A) implies that
\[\LL\simeq \perf(\AA(T^*S^1,D)).\] This category has an
interpretation on the B-side. Consider the singular variety given by
\[ \YY_0= \Spec R[u,v]/(uv-t_0\cdots t_n) \subset \A^{n+3} \]
This is a toric singularity. Indeed, consider the vector space
$V= \A^{2(n+1)}$ generated by the entries of the $2$-by-$(n+1)$
matrix

\begin{center}
\begin{tikzpicture}
  \matrix (m) [matrix of math nodes, nodes in empty cells,
               nodes={minimum width=1.5em, minimum height=2ex, anchor=center},
               column sep=-\pgflinewidth, row sep=-\pgflinewidth, left delimiter={(}, right delimiter={)}
               ]
  {
    x_{0} & x_{1} & \cdots & x_{n} \\
    y_{0} & y_{1} & \cdots & y_{n} \\
  };
\end{tikzpicture}
\end{center}
and consider the action of the torus $T= \mathbb{G}_m^n$
whose $i^{th}$ component acts as follows:
\[ \lambda : \left( {\begin{array}{cccccc} x_{0} & \ldots &
        x_{i-1} & x_{i} & \ldots & x_{n} \\ y_{0} & \ldots &
        y_{i-1} & y_{i} & \ldots & y_{n} \\ \end{array} }
  \right) \to \left( {\begin{array}{cccccc} x_{0} & \ldots &
        \lambda x_{i-1} & \lambda^{-1} x_{i} & \ldots & x_{n} \\
        y_{0} & \ldots & \lambda^{-1} y_{i-1} & \lambda y_{i} &
        \ldots & y_{n} \\ \end{array} } \right) \] Then $\YY_0$
can be identified with the affine GIT quotient $V \sslash T$,
where we can see that $t_i = x_{i}y_{i}$,
$u= x_{0}x_{1}\ldots x_{n}$ and $v= y_{0} y_{1} \ldots,
y_{n}$. The generic GIT quotients $V \sslash_\theta T$ provide
toric crepant resolutions of $\YY_0$. These correspond to
triangulations of $[0,1] \times \Delta_{n}$ where
\(\Delta_{n}\) denotes the \(n\)-simplex. All of these are
(non-canonically) isomorphic to a toric Calabi-Yau variety,
which we denote by $\YY$. These toric Calabi-Yau varieties are
well-known (\cite{DS}, \cite{LPsym2}). We have a map
$\YY \to \Spec R$ given by projection to $(t_0,\ldots t_n)$. The
fiber of this map over 0 is a nodal curve given by a chain of
$\P^1$'s together with two $\A^1$'s attached at the two ends,
and the total space $\YY$ is the versal deformation of this
nodal curve.

There is a tilting bundle $\VV$ on $\YY$ constructed by Van den
Bergh \cite{VdB}; we review this construction in Section
\ref{sct:b_side}. In \ref{cor:universal_calculation}, we will
see that \(\OP{End}_{\YY}(\VV)\) is precisely our algebra
\(\AA(T^*S^1,D)\) and since \(\YY\) is smooth, this means that
\[ \mathcal{L} \simeq D^b(\Coh(\YY)) \] which can be regarded as
a relative version of homological mirror symmetry for $\mathcal{Y}$ (see
also Remark \ref{rem13}).

The braid group action on $D^b (\Coh(\YY))$ is constructed by
Donovan-Segal \cite{DS} by the variation of GIT method, and
previously by Bezrukavnikov-Riche \cite{BR} via Springer
theory. Under the mirror symmetry equivalence discussed above
their action on the $B$-side almost certainly corresponds to our
braid group action on the $A $-side given by Theorem
\ref{thm:main}(B) but we do not check the details here.

\paragraph{Base change.}\label{pg:base_change} We get further
results by working over an \(R\)-algebra \(S\). Let
$\YY_{S,0} = \Spec(\OO_{Y_0}\otimes_R S)$.  Let \(\YY_S\) be the
fibre product:

\begin{center}
  \begin{tikzpicture}
    \node (A) at (0,1) {\(\YY_S\)};
    \node (B) at (2,1) {\(\YY\)};
    \node (C) at (0,0) {\(\YY _{S,0}\)};
    \node (D) at (2,0) {\(\YY_0\)};
    \draw[->] (A) -- (B) node [midway,above] {\(j\)};
    \draw[->] (A) -- (C);
    \draw[->] (C) -- (D);
    \draw[->] (B) -- (D);
    
  \end{tikzpicture}
\end{center}

In \ref{cor:base_change}, we will show that the pullback
\(j^*\VV\) is still a tilting object with
\[\End(j^*\VV)\cong \AA(T^*S^1,D)\otimes_R S.\] The variety
\(\YY_S\) is a partial resolution of \(\YY_{S,0}\), and Theorem
\ref{thm:main}(B) now yields an action of \(\Gamma(T^*S^1,D)\)
by autoequivalences on \(\perf(\YY_S)\). If \(\YY_S\) is itself
smooth, this category is quasi-equivalent to
\(D^b(\Coh(\YY_S))\).

\paragraph{Example.}\label{exm:a_m} If we take \(S=k[t]\)
considered as an \(R\)-module via the homomorphism
\(t_i\mapsto t\) then
\(\YY_{S,0}=\Spec\left(k[u,v,t]/(uv-t^{n+1})\right)\) is the
\(A_n\) surface singularity and \(\YY_S\) is its minimal
resolution, so we get a \(\Gamma(T^*S^1,D)\) action on
\(D^b(\Coh(\YY_S))\). This is one of the examples where we get a
bigger group action: any compactly-supported graded
symplectomorphism of \(T^*S^1\) fixing \(D\) {\em setwise} acts
as an autoequivalence of \(\LL\). This yields an action of the
annular (extended) braid group by autoequivalences. In this
example, an action of the (usual) braid group was known to
Seidel and Thomas \cite{SeidelThomas} and an extended braid
group action was constructed by Gadbled, Thiel and Wagner in
\cite{GadbledThielWagner}.

\paragraph{Example.}\label{exm:can} Let \(f(x,y)\) be a
polynomial whose lowest order term has degree \(n+1\) and
consider the compound \(A_n\) singularity
\(\{uv=f(x,y)\}\subset\C^4\). If \(f\) factors as
\(f_0\cdots f_{n}\) with each curve \(\{f_i(x,y)=0\}\) smooth
then the singularity admits a small resolution. This resolution
has the form \(\YY_S\) where \(S = k[x,y]\) is considered as an
\(R\)-algebra via the homomorphism \(t_i\mapsto f_i(x,y)\). The
algebra \(\AA(T^*S^1,D)\otimes_R S\) is called a {\em
  noncommutative crepant resolution} (NCCR) of this singularity:
it is a noncommutative algebra whose derived category is
equivalent to the derived category of the resolution.

Theorem \ref{thm:main}(B) yields an action of
\(\Gamma(T^*S^1,D)\) on \(D^b(\Coh(\YY_S))\). This can be
enhanced to the bigger group of symplectomorphisms: let \(\psi\)
be a symplectomorphism of \(T^*S^1\) fixing \(D\) setwise and
let \(\sigma\) be the permutation \(\psi(z_i)=z_{\sigma(i)}\);
we get an autoequivalence from \(\psi\) if \(f_{\sigma(i)}=f_i\)
for all \(i\). Autoequivalences of \(D^b(\Coh(\YY_S))\) called
``mutation functors'' were constructed by Iyama and Wemyss
\cite{IW18} using flops along the exceptional curves.

\pg \label{rem13} These examples show that, although this Fukaya
category leaves much of the heavy-lifting to the module category
of the coefficient ring, it does readily give geometric insights
which are nontrivial on the \(B\)-side. The relative Fukaya
category $\WW(T^*S^1, D)$ is appealing because working with
Fukaya categories of surfaces reduces to combinatorial
algebra. However, in view of {\cite[Conjecture E]{LS}}, it is
possible to relate the relative Fukaya category $\WW(T^*S^1, D)$
to an appropriate subcategory of an absolute Fukaya category of
a higher dimensional symplectic manifold $X$. See {\cite[Example
  2.5]{LS}} for a detailed exposition of the case $D=\{1\}$.

\paragraph{Derived contraction algebra.} The derived contraction
algebra is a DG-algebra associated to a small resolution
$\mathcal{Y} \to \mathcal{Y}_0$ that prorepresents derived
deformations of the irreducible components of the reduced
exceptional fiber of the contraction. Concretely, it is a
non-positively graded DG-algebra whose zeroth cohomology
recovers the contraction algebra of Donovan and Wemyss
\cite{DW16}. See the papers by Hua--Toda \cite{HuaToda}, Hua
\cite{Hua}, Hua--Keller \cite{HuaKeller}, and Booth \cite{Booth}
for more background. The derived contraction algebra is obtained
by localising a noncommutative resolution away from an
idempotent. From the Fukaya-categorical description of the
noncommutative resolution in the \(cA_n\) case from
\ref{exm:can}, we can give a geometric interpretation of this
localisation: the derived contraction algebra can be described
using the relative Fukaya category of the punctured disc
\((T^*S^1\setminus L_0,D)\). We discuss this in Section
\ref{sct:dca}.

\paragraph{Acknowledgements.}

JE is supported by EPSRC grant EP/W015749/1. YL is partially
supported by the Royal Society URF\textbackslash R\textbackslash
180024 and EPSRC grant EP/W015889/1. We would like to thank
Michael Wemyss for enlightening discussions which led to a much
cleaner approach, and Matt Booth, Gustavo Jasso, Daniil Mamaev and Richard Thomas for helpful conversations.

\section{The Floer cohomology algebra}
\label{sct:hf_comp}

\begin{figure}[htb]
  \begin{center}
    \begin{tikzpicture}
      \node (A) at (0:2) {\(\bullet\)};
      \node (B) at (60:2) {\(\bullet\)};
      \node (C) at (120:2) {\(\bullet\)};
      \node (D) at (180:2) {};
      \node at (180:2) {\(\cdots\)};
      \node (E) at (240:2) {\(\bullet\)};
      \node (F) at (300:2) {\(\bullet\)};
      \draw[->] (A) to[bend right] (B);
      \draw[->] (B) to[bend right] (C);
      \draw[->] (E) to[bend right] (F);
      \draw[->] (F) to[bend right] (A);
      \draw[->] (C) to[bend right] (D);
      \draw[->] (B) to[bend right] (A);
      \draw[->] (C) to[bend right] (B);
      \draw[->] (F) to[bend right] (E);
      \draw[->] (D) to[bend right] (C);
      \draw[->] (A) to[bend right] (F);
      \draw[->] (D) to[bend right] (E);
      \draw[->] (E) to[bend right] (D);
      \node at (30:2.5) {\(a_{0}\)};
      \node at (90:2.5) {\(a_{1}\)};
      \node at (-90:2.5) {\(a_{n-1}\)};
      \node at (-30:2.55) {\(a_{n}\)};
      \node at (30:1) {\(b_{0}\)};
      \node at (90:1) {\(b_{1}\)};
      \node at (330:1) {\(b_{n}\)};
      \node at (270:0.9) {\(b_{n-1}\)};
      \node at (0:2.4) {\(L_0\)};
      \node at (60:2.4) {\(L_1\)};
      \node at (120:2.4) {\(L_2\)};
      \node at (240:2.4) {\(L_{n-1}\)};
      \node at (300:2.4) {\(L_n\)};   
    \end{tikzpicture}
  \end{center}
  \caption{The quiver \(Q_{n+1}\).}
  \label{fig:quiver}
\end{figure}
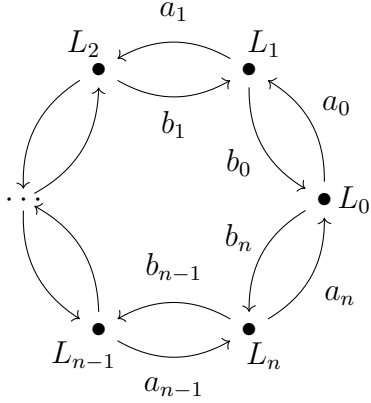

\paragraph{Definition of \(\AA(T^*S^1,D)\).}\label{dfn:alg}
Let \(Q_{n+1}\) be the quiver in Figure \ref{fig:quiver} with
vertices \(L_0,\ldots,L_n\) and arrows\footnote{Indices are
  taken to belong to the cyclic group \(\Z/(n+1)\).}
\(a_{i}\colon L_{i-1}\to L_{i}\),
\(b_{i}\colon L_{i}\to L_{i-1}\).

Recall that \(R=k[t_0,\ldots,t_n]\). Consider the path algebra
\(RQ_{n+1}\) of \(Q_{n+1}\) with coefficients in the ring \(R\); that is
elements of \(RQ_{n+1}\) are \(R\)-linear combinations of paths in
\(Q_{n+1}\) and multiplication is given by concatenate-or-die.  We
write \(e_i\) for the idempotent corresponding to the constant
(lazy) path at the vertex \(L_i\). Let \(I_R\subset RQ_{n+1}\) be
the ideal of \(RQ_{n+1}\) generated by
\[a_{i}b_{i} - t_ie_{i+1},\quad b_{i}a_{i} - t_ie_{i},\quad
  i=0,\ldots,n.\] Write \(\AA(T^*S^1,D)\) for the algebra
\(RQ_{n+1}/I_R\), considered as an \(A_\infty\)-algebra concentrated
in degree zero with no differential or higher operations.

Theorem \ref{thm:main}(A) follows immediately from the next
proposition.

\paragraph{Proposition.}\label{prop:hf_comp} {\it The
  \(A_\infty\)-algebra \(\bigoplus_{i,j=0}^nCF(L_i,L_j)\) is
  quasi-equivalent to \(\AA(T^*S^1,D)\). Note that, in this
  proof, we write \(CF\) to mean \(\hom_{\WW(T^*S^1,D)}\).}
\begin{proof}
  We will use the model of the Fukaya category from
  \cite{HKK}. The arrows labelled \(a\) and \(b\) in Figure
  \ref{fig:quiver} represent the Reeb chords with the same names
  in Figure \ref{fig:sigma}, considered as wrapped intersection
  points \(a_{i}\in CF^0(L_{i},L_{i+1})\),
  \(b_{i}\in CF^0(L_{i+1},L_{i})\). All Reeb chords (called
  ``boundary paths'' in \cite{HKK}) can be obtained by
  concatenating these, and therefore the \(R\)-module
  \(CF(L_i,L_j)\) has as a basis the set of all paths from
  \(L_i\) to \(L_j\) in \(Q_{n+1}\). Here, we include the
  constant path \(e_i\) at \(L_i\), thought of as the identity
  element of \(CF(L_i,L_i)\).

  Since all of these chords are concatenations of chords of
  degree zero, everything is in degree zero, which implies that
  the only nontrivial \(\mu_k\)-operation on
  \(\bigoplus_{i,j}CF(L_i,L_j)\) is \(\mu_2\): the differential
  and higher products all vanish. To compute \(\mu_2\), aside
  from concatenation of chords, we need to count polygons. The
  arcs \(L_i\) cut \(\Sigma\) into \(n+1\) quadrilaterals
  \(D_0,\ldots,D_n\), where we write \(D_i\) for the
  quadrilateral containing the point \(z_i\). Using the
  formula\footnote{The authors of \cite{HKK} state this formula
    for \(\mu_k\) with \(k\geq 3\) only because they do not have
    any quadrilaterals like \(D_i\) in \cite{HKK}.}
  {\cite[Eq. 3.18]{HKK}} and keeping track of our
  additional weighting from the marked points, we see that:
  \[\mu_2(a_{i},b_{i})=t_ie_{i+1}\qquad
    \mu_2(b_{i},a_{i})=t_ie_i\] for all \(i\), where these
  contributions come from \(D_i\). Any other contributions to
  \(\mu_2\) would need to come from quadrilaterals, and any
  quadrilateral can be decomposed as a union of \(D_i\)s, so any
  other \(\mu_2\) product can be deduced from these.
\end{proof}

\section{Autoequivalences}\label{sct:group_action}

\paragraph{Group action.} \label{pg:group_action} Let
\(R=k[t_0,\ldots,t_n]\). Given a permutation \(\sigma\) of
\(\{0,1,\ldots,n\}\), let \(R_\sigma\) denote the \(R\)-module
whose underlying vector space is \(R\) but \(t_i\) acts as
multiplication by \(t_{\sigma(i)}\). Consider the triangulated
\(A_\infty\)-category
\[\WW(T^*S^1,D)\rtimes S_{n+1} := \coprod_{\sigma\in S_{n+1}} \WW(T^*S^1,D)\otimes_R R_\sigma\]
where the morphism spaces between different components are
zero. Given a graded symplectomorphism
\(\psi\colon T^*S^1\to T^*S^1\) satisfying \(\psi(D)=D\), we get
a permutation \(\sigma\in S_{n+1}\) defined by
\(\psi(z_i)=z_{\sigma(i)}\). This induces an autoequivalence
\[\WW(T^*S^1,D)\rtimes S_{n+1}\to \WW(T^*S^1,D)\rtimes S_{n+1}\] sending
\(\WW(T^*S^1,D)\otimes_R R_{\tau}\) to
\(\WW(T^*S^1,D)\otimes_R R_{\sigma\tau}\). In particular, this
gives an action of the pure annular braid group by
autoequivalences on \(\WW(T^*S^1,D)\).

\paragraph{Theorem.} {\em Let \(\LL_\sigma\) denote the
  subcategory of \(\WW(T^*S^1,D)\otimes_R R_\sigma\) generated
  by the arcs \(L_0,\ldots,L_n\). Then the autoequivalences from
  \ref{pg:group_action} preserve
  \(\coprod_{\sigma\in S_{n+1}}\LL_\sigma\).}

We now begin the proof of this theorem, which will conclude in
\ref{pg:end_of_proof}. We will focus on the case \(n\geq 2\)
because it can be handled uniformly: for small \(n\) the
arguments are similar but the pictures are slightly different
because \(L_1=L_n\) or \(L_0=L_1=L_n\). Throughout the argument
we will ignore signs and orientations of moduli spaces. The
reason we can get away with this is explained in Remark
\ref{rmk:signs}.
 
\pg\label{pg:symplectos} We define some compactly-supported
symplectomorphisms of \(T^*S^1\) fixing \(D\) setwise. First,
let \((p,q)\) be coordinates with \(p\in\R\) and
\(q\in S^1=\R/2\pi\Z\), and define the symplectomorphisms
\[\rho(p,q)=(p,q+f(p)),\qquad \delta(p,q)=(p,q+g(p))\]
where \(f,g\colon\R\to\R\) are the functions shown in Figure
\ref{fig:functions}.

\begin{figure}[htb]
  \begin{center}
    \begin{tikzpicture}
      \draw[->] (-2,0) -- (2,0) node [right] {\(p\)};
      \draw[->] (0,0) -- (0,2) node [above] {\(f\)};
      \draw[dotted] (0,2) -- (-2,2) node [left] {\(2\pi\)};
      \draw[dotted] (0,1) -- (-2,1) node [left] {\(\pi\)};
      \draw[dotted] (-0.75,1) -- (-0.75,0) node [below] {\(p_0\)};
      \draw[thick] (-2,0) -- (-1.5,0) to[out=0,in=180] (0,2) -- (2,2);
      \begin{scope}[shift={(6,0)}]
      \draw[->] (-2,0) -- (2,0) node [right] {\(p\)};
      \draw[->] (0,0) -- (0,2) node [above] {\(g\)};
      \draw (-2,0) -- (2,0);
      \draw (0,0) -- (0,2);
      \draw[dotted] (-2,1) -- (2,1) node [right] {\(2\pi/(n+1)\)};
      \draw[thick] (-2,0) -- (-1.5,0) to[out=0,in=180] (0,1) to[out=0,in=180] (1.5,0) -- (2,0);
      \end{scope}
      
    \end{tikzpicture}
    \caption{The functions \(f\) and \(g\) used in the
      definitions of the symplectomorphisms \(\rho\) and
      \(\delta\) in \ref{pg:symplectos}.}
    \label{fig:functions}
  \end{center}
\end{figure}

The symplectomorphism \(\rho\) fixes the two noncompact ends and
rotates the points in \(D\) by an angle \(2\pi/(n+1)\); the
symplectomorphism \(\delta\) is a Dehn twist along a loop
\(\{p_0\}\times S^1\) with \(p_0<0\). Next, let
\(\psi_i\colon T^*S^1\to T^*S^1\) denote the half-twist around
the arc connecting \(z_{i-1}\) to \(z_{i}\) (indices taken
modulo \(n+1\)). The mapping classes
\(\psi_0,\ldots,\psi_n,\rho,\delta\) generate the graded
symplectic mapping class group: see\footnote{Gadbled, Thiel and
  Wagner treat one of the two noncompact ends as a puncture, so
  do not need \(\delta\).} {\cite[Section
  1]{GadbledThielWagner}}. The symplectomorphism \(\delta\) acts
trivially on our Lagrangians as objects of the wrapped category:
\(\delta\) is part of the wrapping that we would do anyway to
compute \(\hom\)-spaces. The symplectomorphism \(\rho\)
cyclically permutes the \(L_i\) (up to Hamiltonian isotopy). So
to prove that \(\Gamma(T^*S^1,D)\) preserves \(\LL\), it
suffices to check that \(\psi_i(L_j)\) is generated by the arcs
\(L_0,\ldots,L_n\) for all \(i,j\). In fact, \(\psi_i(L_j)=L_j\)
unless \(i=j\), so we just need to study
\(\psi_i(L_i)\). Moreover, by cyclic symmetry of \((T^*S^1,D)\)
we can assume that \(i=0\).

\begin{figure}[htb]
  \begin{center}
    \begin{tikzpicture}
      \filldraw[fill=gray, opacity=0.5,draw=none] (120:3) -- (120:1) arc [radius = 1,start angle=120,end angle=0] -- (2,0) to[out=135,in=-90] (80:2) to[out=90,in=180] (60:3) arc [radius = 3,start angle=60,end angle =120];
      \filldraw[fill=gray, opacity=0.5,draw=none] (-120:1) -- (-120:3) arc [radius = 3,start angle=-120,end angle=0] -- (2,0) to[out=-45,in=0] (-90:2.5) to[out=180,in=-90] (-60:1) arc [radius = 1,start angle=-60,end angle =-120];
      \node at (100:2) {\(A\)};
      \node at (-105:2) {\(B\)};
      \draw (0,0) circle [radius = 3];
      \draw (0,0) circle [radius = 1];
      \draw (0:1) -- (0:3);
      \draw (120:1) -- (120:3);
      \draw (-120:1) -- (-120:3);
      \node at (-120:3) [below left] {\(L_n\)};
      \node at (120:3) [above left] {\(L_1\)};
      \node at (3,0) [right] {\(L_0\)};
      \node (z2) at (60:2) {\(\star\)};
      \node (z1) at (-60:2) {\(\star\)};
      \node at (z1) [above] {\(z_n\)};
      \node at (z2) [below] {\(z_0\)};
      \draw (60:3) to[out=180,in=90] (80:2) to[out=-90,in=135] (0:2) to[out=-45,in=0] (-90:2.5) to[out=180,in=-90] (-60:1);
      \node at (175:2) {\(\cdot\)};
      \node at (180:2) {\(\cdot\)};
      \node at (185:2) {\(\cdot\)};
      \node at (180:3) {\(\circ\)};
      \node at (180:1) {\(\circ\)};
      \draw [-{Stealth[length=3mm]}] (0:3) arc [radius = 3,start angle = 0, end angle = 30] node [right] {\(\alpha\)};
      \draw [-{Stealth[length=3mm]}] (60:3) arc [radius = 3,start angle = 60, end angle = 90] node [above right] {\(\alpha'\)};
      \draw [-{Stealth[length=3mm]}] (-120:3) arc [radius = 3,start angle = -120, end angle = -60] node [below right] {\(a_n\)};
      \draw [-{Stealth[length=3mm]}] (120:1) arc [radius = 1,start angle = 120, end angle = 60] node [below] {\(b_0\)};
      \draw [-{Stealth[length=3mm]}] (0:1) arc [radius = 1,start angle = 0, end angle = -30];
      \node at (-30:0.8) {\(\beta\)};
      \draw [-{Stealth[length=3mm]}] (-60:1) arc [radius = 1,start angle = -60, end angle = -90];
      \node at (-90:0.8) {\(\beta'\)};
      \node at (2,0) [above] {\(p\)};
      \node at (2,0) {\(\bullet\)};
      \node at (62:3.1) [right] {\(\psi_0(L_0)\)};
    \end{tikzpicture}
    \caption{The half-twisted arc \(\psi_0(L_0)\), perturbed
      slightly along the Reeb flow to separate it from
      \(L_0\). We have added two stops on the boundary for
      convenience; these are labelled \(\circ\). We have also
      labelled the Reeb orbits connecting the Lagrangian
      arcs. Note that \(a_0=\alpha'\alpha\) and
      \(b_n=\beta'\beta\). The point \(p\) (marked with a
      \(\bullet\)) is an intersection point of \(L_0\) with
      \(\psi_0(L_0)\). Two important polygonal regions \(A\) and
      \(B\) are shaded.}
    \label{fig:half_twist}
  \end{center}
\end{figure}
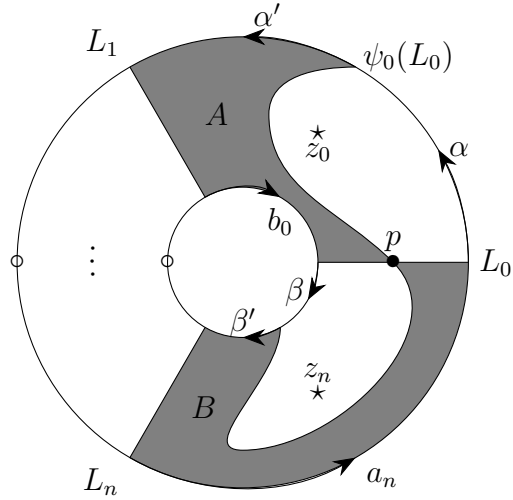

\pg The half-twisted arc \(\psi_0(L_0)\) is shown in Figure
\ref{fig:half_twist}. To localise the calculation near the
diagram, we will insert a stop (in the sense of Sylvan
\cite{Sylvan}) on each of the two boundary components and work
first in the partially wrapped Fukaya category. We will write
down a twisted complex \(\L'\) built out of \(L_n\), \(L_0\)
and \(L_1\) and a quasi-isomorphism
\(q\in CF(\L',\psi_0(L_0))\). If we then apply Sylvan's stop
removal functor to this twisted complex, we obtain a twisted
complex \(\L\) in \(\WW(T^*S^1,D)\) which is quasi-isomorphic
to \(\psi_0(L_0)\).

\pg  The advantage of inserting stops is that the partially wrapped
Floer cohomology is easy to read off from Figure
\ref{fig:half_twist}:
\begin{align*}
  CF(\psi_0(L_0),L_0)&= R\cdot p,& CF(L_0,\psi_0(L_0))&=R\cdot
                                                        p\,\oplus\,
                                                        R\cdot\alpha\,\oplus\,
                                                        R\cdot\beta,\\
  CF(\psi_0(L_0),L_1)&=R\cdot\alpha',&CF(L_1,\psi_0(L_0))&=R\cdot
                                                           (\beta b_0),\\
  CF(\psi_0(L_0),L_n) &=R\cdot\beta',&CF(L_n,\psi_0(L_0))&=R\cdot(\alpha a_n).
\end{align*}
All of these morphisms are in degree zero except for \(p\)
which is in degree \(1\).

\pg  Consider the twisted complex
\[\L' \coloneqq \left(L_1\oplus
    L_n\xrightarrow{(b_0,a_n)} L_0\right)\] and the morphisms
\(q_1\colon \L'\to \psi_0(L_0)\) and \(q_2\colon
\psi_0(L_0)\to \L'\) defined by\footnote{We will write twisted
  complexes horizontally and morphisms between them vertically.}

\begin{center}
  \begin{tikzpicture}
    \node at (-1,-1) {\(q_1\colon\)};
    \node (a) at (0,0) {\(L_1\oplus L_n\)};
    \node (b) at (3,0) {\(L_0\)};
    \draw[->] (a) -- (b) node [midway,above] {\((b_0,a_n)\)};
    \node (c) at (0,-2) {\(\psi_0(L_0)\)};
    \draw[->] (b) -- (c) node [midway, below right] {\(p\)};
    \begin{scope}[shift={(6,0)}]
      \node at (-1,-1) {\(q_2\colon\)};
      \node (a2) at (0,-2) {\(L_1\oplus L_n\)};
      \node (b2) at (3,-2) {\(L_0\)};
      \draw[->] (a2) -- (b2) node [midway,below] {\((b_0,a_n)\)};
      \node (c2) at (0,0) {\(\psi_0(L_0)\)};
      \draw[->] (c2) -- (a2) node [midway, right] {\((\alpha',\beta')\)};
    \end{scope}
  \end{tikzpicture}
\end{center}

We need to show that \(\mu^{Tw}_2(q_1,q_2)\) and
\(\mu^{Tw}_2(q_2,q_1)\) are equal to the identity elements of
\(CF(\psi_0(L_0),\psi_0(L_0))\) and \(CF(\L',\L')\)
respectively (we are using Seidel's convention for
composition, right-to-left). We compute \(\mu^{Tw}_2\) by
stacking the morphisms and then taking all possible paths
through the resulting diagram, composing wherever possible.

\pg  To calculate \(\mu^{Tw}_2(q_2,q_1)\), we have the following
diagram:

\begin{center}
  \begin{tikzpicture}
    \node (d) at (0,-2) {\(\psi_0(L_0)\)};
    \node (a) at (0,0) {\(L_1\oplus L_n\)};
    \node (b) at (3,0) {\(L_0\)};
    \draw[->] (a) -- (b) node [midway,above] {\((b_0,a_n)\)};
    \node (a2) at (0,-4) {\(L_1\oplus L_n\)};
    \node (b2) at (3,-4) {\(L_0\)};
    \draw[->] (a2) -- (b2) node [midway,below] {\((b_0,a_n)\)};
    \draw[->] (b) -- (c) node [midway, below right] {\(p\)};
    \draw[->] (d) -- (a2) node [midway, left] {\((\alpha',\beta')\)};
  \end{tikzpicture}
\end{center}

There are several routes through the diagram connecting the
top row to the bottom. There are two paths that involve three
morphisms:

\begin{center}
  \begin{tikzpicture}
    \node (a) at (0,0) {\(L_1\oplus L_n\)};
    \node (b) at (3,0) {\(L_0\)};
    \draw[->] (a) -- (b) node [midway,above] {\((b_0,a_n)\)};
    \node (a2) at (0,-2) {\(L_1\oplus L_n\)};
    \node (b2) at (3,-2) {\(L_0\)};
    \draw[->] (a2) -- (b2) node [midway,below] {\((b_0,a_n)\)};
    \draw[->] (a) -- (a2);
    \node at (0,-1) [left] {\(\begin{pmatrix}
        \mu_3(\alpha',p,b_0) & \mu_3(\alpha',p,a_n) \\
        \mu_3(\beta',p,b_0) & \mu_3(\beta',p,a_n) \end{pmatrix}\)};
    \draw[->] (b) -- (b2);
    \node at (3,-1) [right] {\(\mu_3(b_0,\alpha',p)+\mu_3(a_n,\beta',p)\)};
  \end{tikzpicture}
\end{center}

There is also a path of length \(2\) connecting \(L_0\) to
\(L_1\oplus L_n\) and one of length \(4\) connecting
\(L_1\oplus L_n\) to \(L_0\). Both of these concatenations
vanish for degree reasons.

\begin{figure}[htb]
  \begin{center}
    \begin{tikzpicture}[scale=1.5]
      \draw (0,0) circle [radius = 3];
      \draw (0,0) circle [radius = 1];
      \node (alphaprime) at (145:2.6) {};
      \node (alpha) at (90:2.3) {};
      \node (beta) at (-90:1.9) {};
      \node (betaprime) at (-120:1.6) {};
      \node (b0) at (1.2,0) {};
      \node (pold) at (2,0.5) {};
      \node (b000) at (-1.0625,-0.273) {};
      \node (annnlabel) at (-2.501184896,0.85057) {};
      \node (annn) at (165:2.9) {};
      \node at (180:3) {\(\circ\)};
      \node at (180:1) {\(\circ\)};
      \node (idb) at (-0.736,-0.817) {};
      \node (ida) at (-2.014,2.086) {};
      \node (pactual) at (-25:1.9) {};
      \begin{scope}[shift={(-1,4.5)}]
        \draw[pattern=north west lines, pattern color=gray,opacity=0.5,draw=none] (-5,-3) rectangle (-4,-2);
        \draw[blue] (-5,-3) -- (-4,-3) node[midway,below] {\(L'''_0\)};
        \draw[purple] (-5,-3) -- (-5,-2) node[midway,left] {\(L'''_1\)};
        \draw[black] (-5,-2) -- (-4,-2) node[midway,above] {\(L'_1\)};
        \draw[red] (-4,-3) -- (-4,-2) node [midway,right] {\(\psi_0(L_0)''\)};
        \node at (-5,-3) [above right] {\footnotesize \(b'''_0\)};
        \node at (-5,-2) [below right] {\footnotesize \(\otimes\)};
        \node at (-4,-3) [above left] {\footnotesize \(p\)};
        \node at (-4,-2) [below left] {\footnotesize \(\alpha'\)};
      \end{scope}
      \begin{scope}[shift={(-1,2.5)}]
        \draw[pattern=north east lines, pattern color=gray,opacity=0.5,draw=none] (-5,-3) rectangle (-4,-2);
        \draw[blue] (-5,-3) -- (-4,-3) node[midway,below] {\(L'''_0\)};
        \draw[purple] (-5,-3) -- (-5,-2) node[midway,left] {\(L'''_n\)};
        \draw[black] (-5,-2) -- (-4,-2) node[midway,above] {\(L'_n\)};
        \draw[red] (-4,-3) -- (-4,-2) node [midway,right] {\(\psi_0(L_0)''\)};
        \node at (-5,-3) [above right] {\footnotesize \(a'''_n\)};
        \node at (-5,-2) [below right] {\footnotesize \(\otimes\)};
        \node at (-4,-3) [above left] {\footnotesize \(p\)};
        \node at (-4,-2) [below left] {\footnotesize \(\beta'\)};
      \end{scope}
      \begin{scope}[shift={(-1,0.5)}]
        \draw[pattern=vertical lines, pattern color=gray,opacity=0.5,draw=none] (-5,-3) rectangle (-4,-2);
        \draw[blue] (-5,-3) -- (-4,-3) node[midway,below] {\(L'''_0\)};
        \draw[red] (-5,-3) -- (-5,-2) node[midway,left] {\(\psi_0(L_0)''\)};
        \draw[black] (-5,-2) -- (-4,-2) node[midway,above] {\(L_0\)};
        \draw[black] (-4,-3) -- (-4,-2) node [midway,right] {\(L'_n\)};
        \node at (-5,-3) [above right] {\footnotesize \(p\)};
        \node at (-5,-2) [below right] {\footnotesize \(\beta\)};
        \node at (-4,-3) [above left] {\footnotesize \(\otimes\)};
        \node at (-4,-2) [below left] {\footnotesize \(a_n\)};
      \end{scope}
      \filldraw[pattern=north west lines, pattern color=gray,opacity=0.5,draw=none] (b000.center) arc[radius=1.1,start angle=-165,end angle=-132.3] to[out=-90,in=150] (-120:1.2) arc[radius=1.2,start angle=-120,end angle=105] to[out=-165,in=-30] (alphaprime.center) to[out=55,in=90] (1.7,0) -- (pactual.center) to[out=-118,in=0] (beta.center) to[out=180,in=-90] (-160:1.8) to[out=90,in=-170] (b000.center);
      \filldraw[pattern=north east lines, pattern color=gray,opacity=0.5,draw=none] (annn.center) to[out=-10,in=-105] (165:2.3) arc [radius=2.3,start angle=165,end angle=60] to[out=-30,in=90] (1.9,0) to[out=-90,in=70] (pactual.center) to[out=-90,in=0] (1,-2.2) to[out=180,in=-30] (betaprime.center) to[out=-135,in=150] (-120:2.8) arc [radius=2.8,start angle=-120,end angle=130] to[out=-130,in=-45] (135:2.9) arc [radius=2.9, start angle=135,end angle=145];
      \filldraw[pattern=vertical lines, pattern color=gray,opacity=0.5,draw=none] (1.9,0) -- (2.8,0) arc[radius=2.8,start angle = 0,end angle=-120] to[out=150,in=-135] (betaprime.center) to[out=-30,in=180] (1,-2.2) to[out=0,in=-87] (pactual.center) to[out=70,in=-90] (1.9,0);
      \draw (135:3) to[out=-45,in=-130] (130:2.8) arc [radius=2.8,start angle=130,end angle=-120] to[out=150,in=-135] (betaprime.center) to[out=45,in=-135] (-145:1);
      \node at (130:3.1) [left] {\(L'_n\)};
      \node at (2.8,0) {\footnotesize \(\bullet\)};
      \node at (2.65,0.2) {\scriptsize \(a_n\)};
      \draw (145:3) to[out=-35,in=-165] (105:1.2) arc [radius=1.2,start angle=105,end angle=-120] to[out=150,in=-135] (-135:1); \node at (142.5:3.1) [left] {\(L'_1\)};
      \draw (1,0) -- (3,0); \node at (3,0) [right] {\(L_0\)};
      \draw[red] (155:3) to[out=-25,in=-120] (150:2.6) to[out=60,in=-145] (145:2.6) to[out=55,in=90] (1.7,0) to[out=-90,in=0] (1,-2.2) to[out=180,in=-30] (betaprime.center) to[out=150,in=-155] (-155:1);
      \node at (152.5:3) [left] {\(\psi_0(L_0)''\)};
      \draw[blue] (165:3) to[out=-10,in=-105] (165:2.3) arc [radius=2.3,start angle=165,end angle=60] to[out=-30,in=90] (1.9,0) to[out=-90,in=0] (beta.center) to[out=180,in=-90] (-160:1.8) to[out=90,in=-170] (-165:1);
      \node (l000) at (-3.5,1) {\(L'''_0\)};
      \draw[blue,->] (l000) -- (165:3);
      \node (z1) at (-60:2.1) {\(\star\)}; \node (z2) at (60:2) {\(\star\)}; \node at (z1) [below] {\(z_n\)};
      \node at (z2) [below] {\(z_0\)};
      \node at (alphaprime.center) {\footnotesize \(\bullet\)};
      \node at (b0.center) {\footnotesize \(\bullet\)};
      \node (id) at (1.9,0) {\footnotesize \(\bm{\otimes}\)};
      \node (alphaactual) at (117:2.3) {\footnotesize \(\bullet\)};
      \node at (alphaactual) [below] {\scriptsize \(\alpha\)};
      \node at (pactual) {\footnotesize \(\bullet\)};
      \node at (pactual) [right] {\scriptsize \(p\)};
      \node at (beta) {\footnotesize \(\bullet\)};
      \node at (beta) [above] {\scriptsize \(\beta\)};
      \node at (1.4,-0.23) {\scriptsize \(b_0\)};
      \node (alphaplabel) at (-2.1297953,1.7212987) {\scriptsize \(\alpha'\)};
      \draw[purple] (170:3) to[out=-10,in=-105] (170:2.9) arc [radius=2.9,start angle=170,end angle=-170] to[out=100,in=180] (-176:1);
      \node (lnnnn) at (-3.5,0.5) {\(L''''_n\)};
      \draw[purple,->] (lnnnn) -- (170:3);
      \draw[purple] (175:3) to[out=-10,in=205] (115:1.1) arc [radius =1.1,start angle= 115,end angle=-165] to[out=105,in=-170] (-173:1);
      \node (l1111) at (-3.5,0) {\(L_1''''\)};
      \draw[purple,->] (l1111) -- (175:3);
      \node at (betaprime.center) {\footnotesize \(\bullet\)};
      \node at (-0.55,-1.2856406) {\scriptsize \(\beta'\)};
      \node at (annn) {\footnotesize \(\bullet\)};
      \node at (annnlabel) {\scriptsize \(a'''_n\)};
      \node at (b000) {\footnotesize \(\bullet\)};
      \node at (b000.south west) [] {\scriptsize \(b'''_0\)};
      \node at (idb.center) {\footnotesize \(\bm{\otimes}\)};
      \node at (ida.center) {\footnotesize \(\bm{\otimes}\)};
    \end{tikzpicture}
    \caption{The choices of partially wrapped Hamiltonian perturbations for the computations in \ref{pg:perturbations}. The intersection points marked \(\bm{\otimes}\) denote the identity elements of the corresponding Floer complex. We show the three holomorphic quadrilaterals which contribute to \(\mu_3(\alpha',p,b'''_0)\), \(\mu_3(\beta',p,a'''_n)\) and \(\mu_3(a_n,\beta',p)\) (all other products vanish with these choices); the quadrilaterals are distinguished by the direction of their hatching.}
    \label{fig:choices}
  \end{center}
\end{figure}

\pg\label{pg:perturbations} Up until this point, we have been
relaxed about choosing Hamiltonian perturbations, but in order
to proceed we must specify which choices of partially wrapped
perturbations have been made. The relevant perturbations are
\(L'_1\), \(L'_n\), \(\psi_0(L_0)''\), \(L'''_0\), \(L''''_1\)
and \(L''''_n\) where each prime indicates that we have wrapped
more; see Figure \ref{fig:choices} for our
specific choices and the relevant intersection points. Note that
we now need to distinguish notationally between
\(b_0\in CF(L_1',L_0)\) and \(b'''_0\in CF(L''''_1,L'''_0)\) and
between \(a_n\in CF(L'_n,L_0)\) and
\(a'''_0\in CF(L''''_n,L'''_0)\). This allows us to read off all
the relevant \(\mu_3\) products contributing to
\(\mu_2^{Tw}(q_2,q_1)\) from quadrilaterals in the picture. The
result is:

\begin{center}
  \begin{tikzpicture}
    \node (a) at (0,0) {\(L_1\oplus L_n\)};
    \node (b) at (3,0) {\(L_0\)};
    \draw[->] (a) -- (b) node [midway,above] {\((b_0,a_n)\)};
    \node (a2) at (0,-2) {\(L_1\oplus L_n\)};
    \node (b2) at (3,-2) {\(L_0\)};
    \draw[->] (a2) -- (b2) node [midway,below] {\((b_0,a_n)\)};
    \draw[->] (a) -- (a2);
    \node at (0,-1) [left] {\(\begin{pmatrix}
        1 & 0 \\
        0 & 1 \end{pmatrix}\)};
    \draw[->] (b) -- (b2);
    \node at (3,-1) [right] {\(1\)};
  \end{tikzpicture}
\end{center}

For example, let us compute \(\mu_3(\alpha',p,b'''_0)\) and see
that it is equal to \(1\). We must think of this \(\mu_3\)
product as a map
\[\mu_3\colon CF(\psi_0(L_0)'',L'_1)\otimes
  CF(L'''_0,\psi_0(L_0)'')\otimes CF(L''''_1,L'''_0)\to
  CF(L''''_1,L'_1).\] In Figure \ref{fig:choices} there is a
unique quadrilateral with vertices at \(\alpha'\), \(p\),
\(b'''_0\), and at the unique intersection point
\(L'''_1\cap L_1\) which represents \(1\in
CF(L'''_1,L_1)\). This shows that
\(\mu_3(\alpha',p,b'''_0)=1\). The other calculations are
similar; note that \(\mu_3(a_n,\beta',p)=1\) and
\(\mu_3(b_0,\alpha',p)=0\) with our choice of perturbations, so
that \(\mu_3(a_n,\beta',p)+\mu_3(b_0,\alpha',p)=1\).

\pg \label{pg:end_of_proof} To calculate
\(\mu^{Tw}_2(q_1,q_2)\), we have the following diagram:

\begin{center}
  \begin{tikzpicture}
    \node (d) at (0,2) {\(\psi_0(L_0)\)};
    \node (a) at (0,0) {\(L_1\oplus L_n\)};
    \node (b) at (3,0) {\(L_0\)};
    \draw[->] (a) -- (b) node [midway,above] {\((b_0,a_n)\)};
    \node (c) at (0,-2) {\(\psi_0(L_0)\)};
    \draw[->] (b) -- (c) node [midway, below right] {\(p\)};
    \draw[->] (d) -- (a) node [midway, left] {\((\alpha',\beta')\)};
  \end{tikzpicture}
\end{center}

There is only one route from the
top row to the bottom, which means that
\begin{align*}
  \mu^{Tw}_2(q_1,q_2)&=\mu_3(p,(b_0,a_n),(\alpha',\beta'))\\
                     &=\mu_3(p,b_0,\alpha') + \mu_3(p,a_n,\beta')
\end{align*} As with the previous calculation, this yields
\(1\in CF(\psi_0(L_0),\psi_0(L_0))\). This shows that \(q_1\) and \(q_2\)
are mutually inverse quasi-isomorphisms, which completes the
proof.\qed

\paragraph{Remark about signs.}\label{rmk:signs} In this proof,
we completely
ignored signs. If we insert all the undetermined signs, the
arguments yield
\[\mu^{Tw}_2(q_1,q_2)=\pm\OP{id}_{\psi_0(L_0)},\qquad\mu^{Tw}_{q_1}(q_2,q_2)
  = (\pm \OP{id}_{L_1})\oplus (\pm \OP{id}_{L_n})\oplus(\pm
  \OP{id}_{L_0}).\] At this point, we pass to cohomology and
consider the morphisms \([q_1]\in HF(\L,\psi_0(L_0))\) and
\([q_2]\in HF(\psi_0(L_0),L_0)\). The morphisms
\[[q_1]\in HF(\L,\psi_0(L_0)),\qquad
  [q_2]\circ[q_1]\circ[q_2]\in HF(\psi_0(L_0),\L)\] are now
mutually inverse because all signs are squared in the composites
\([q_1]\circ[q_2]\circ[q_1]\circ[q_2]\) and
\([q_2]\circ[q_1]\circ[q_2]\circ[q_1]\).

\section{B-side}\label{sct:b_side}

\paragraph{Setup.} As in the introduction, let
\(R=k[t_0,\ldots,t_n]\), let $\YY_0 = \Spec R[u,v]/(uv-t_0\cdots
t_n))$, and let \(f\colon \YY_0 \to \A^{n+1}\) be the morphism given by
$(t_0,\ldots, t_n)$. This morphism \(f\) is the versal deformation of
the \(A_{n}\) curve singularity. We have a toric crepant resolution
$\pi\colon \mathcal{Y} \to \mathcal{Y}_0$ given by a triangulation of
$[0,1] \times \Delta_{n}$.

\paragraph{The Van den Bergh tilting bundle.} We now describe a
tilting bundle on \(\YY\), making explicit the construction of
Van den Bergh {\cite[Propositions 3.2.5, 3.2.10]{VdB}} in this
example. Recall from \ref{pg:mirror_interp} that \(\YY\) is the
GIT quotient \(V\sslash_\theta T\), where \(V\) is the space of
\(2\)-by-\((n+1)\) matrices
\[\begin{pmatrix} x_{0} & \ldots & x_{i} & x_{i+1} &
    \ldots & x_{n} \\ y_{0} & \ldots & y_{i} & y_{i+1} &
    \ldots & y_{n} \\ \end{pmatrix}\] and the torus
\(T=\Gm^{n}\) acts as
\[\begin{pmatrix} \lambda_1x_{0} & \ldots &
    \lambda_{i}^{-1}\lambda_{i+1}x_{i} &
    \lambda_{i+1}^{-1}\lambda_{i+2}x_{i+1} & \ldots &
    \lambda_{n}^{-1}x_{n} \\ \lambda_1^{-1}y_{0} & \ldots
    & \lambda_{i}\lambda_{i+1}^{-1}y_{i} &
    \lambda_{i+1}\lambda_{i+2}^{-1}y_{i+1} & \ldots &
    \lambda_{n}y_{n} \\ \end{pmatrix}\] and \(\theta\) is
the character \(\theta(\lambda_1,\ldots,\lambda_{n}) =
\lambda_1\cdots\lambda_{n}\) of \(T\).

Given another character \(\chi\colon T\to\C^*\), we get a line
bundle \((V\times\C)\sslash_\theta T\) over \(\YY\), where \(T\)
acts with weight \(\chi\) on \(\C\). Let \(\MM_i\) be the line
bundle corresponding to the character
\(\chi_i(\lambda_1,\ldots,\lambda_{n})=\lambda_i\). The
sections of \(\MM_i\) are in bijection with the polynomials in
the variables \( x_{i}, y_{i} \) which have weight \(\chi\) under the
action of \(T\). For example, \(x_{0}\) is a section of
\(\MM_1\) and \(y_n\) is a section of \(\MM_n\).

\paragraph{Lemma.}\label{lma:sections_polynomials} {\em The
  sections of \(\OO_{\YY}\) form a ring isomorphic to
  \(R[u,v]/(uv-t_0\cdots t_n)\). The sections of \(\MM_i\) form
  a module over this ring which is generated
  by \(\sigma_i \coloneqq x_{0}\cdots x_{i-1}\) and
  \(\tau_i \coloneqq y_{i}\cdots y_{n}\).}

Note that since \(\pi_*\OO_\YY=\OO_{\YY_0}\) we can think of
\(H^0(\MM_i)\) as an \(\OO_\YY\)-module or an
\(\OO_{\YY_0}\)-module. It is isomorphic to the
\(R[u,v]/(uv-t_0\cdots t_n)\)-module \((u,t_0\cdots t_{i-1})\) by
identifying \(\sigma_i\) with \(u\) and \(\tau_i\) with
\(t_0\cdots t_{i-1}\).

\begin{proof}
  Consider the monomial
  \(x_{0}^{c_{0}}\cdots x_{n}^{c_{n}}y_{0}^{d_{0}}\cdots
  y_{n}^{d_{n}}\). The condition that this defines a section of
  \(\OO_\YY\) is that \(c_{i}+d_{i+1}-c_{i+1}-d_{i}=0\) for all
  \(i=0,\ldots,n-1\). This implies that
  \(c_{0}-d_{0}=\cdots=c_{n}-d_{n}\). If this common value is
  positive then the monomial can be written as
  \[t_0^{d_{0}}\cdots t_n^{d_{n}}u^{c_{0}-d_{0}}\]
  otherwise it can be written as
  \[t_0^{c_{0}}\cdots t_n^{c_{n}}v^{d_{0}-c_{0}}\]
  where we are defining
  \[u=x_{0}\cdots x_{n},\qquad v=y_{0}\cdots
    y_{n},\qquad t_i=x_{i}y_{i}\] as in
  \ref{pg:mirror_interp}. The argument for the sections of
  \(\MM_i\) is similar except one is left with an additional
  factor of \(x_{0}\cdots x_{i-1}\) or
  \(y_{i+1}\cdots y_{n}\) depending on whether
  \(c_{i}>d_{i}\) or \(d_{i+1}>c_{i}\).
\end{proof}

\paragraph{Lemma.} {\em Let \(\MM=\bigoplus_{i=1}^{n}
\MM_i\). Consider the \(n-1\) sections
\begin{align*}
  s_1&=(\sigma_1,\tau_2,0,\ldots,0)\\
  s_2&=(0,\sigma_2,\tau_3,0,\ldots,0)\\
  &\ \ \ \vdots\\
  s_{n-1}&=(0,\cdots,0,\sigma_{n-1},\tau_{n})
\end{align*}
These sections are everywhere linearly independent, and hence
span a copy of the trivial bundle of rank \(n-1\) inside
\(\MM\).}

\begin{proof}
  At each point of \(\YY\), the wedge product
  \(s_1\wedge s_2\wedge \cdots\wedge s_{n-1}\) has components
  \begin{gather*}
    \tau_2\cdots\tau_{n},\\
    \sigma_1\tau_3\cdots\tau_{n},\\
    \sigma_1\sigma_2\tau_4\cdots\tau_{n},\\
    \vdots\\
    \sigma_1\cdots\sigma_{n-1}.
  \end{gather*}
  If the sections are linearly dependent somewhere then all of
  these components vanish at that point. Let \(j\) be minimal
  such that \(\sigma_j=0\); note that this implies
  \(x_{j}=0\). Since
  \(\sigma_1\cdots\sigma_{j-1}\tau_{j+1}\cdots\tau_{n}=0\) we
  deduce that some \(\tau_k=0\) for \(k>j\), and for the maximal
  such $k$ we have that \(y_{k}=0\). But, as can be easily verified using the Hilbert-Mumford criterion (cf. \cite{thomas}), the unstable locus for the linearization $\theta$ is the union of the subvarieties
  \(\{x_{j} = y_{k} = 0\}\) for \(0\leq j<k\leq n\), so on the
  GIT quotient \(\YY\) there are no points where these sections
  vanish simultaneously.
\end{proof}

\paragraph{Corollary.} {\em Let \(\LL\) be the quotient of
  \(\MM\) by the trivial subbundle spanned by these
  sections. Then \(\LL\) is an ample line bundle on \(\YY\) and
  \(\VV\coloneqq\OO_\YY\oplus\MM\) is a tilting bundle.}

\begin{proof}
  The quotient is a line bundle and is therefore determined by
  its first Chern class, which is in turn determined by its
  restriction to the curve
  \(\{t_0=\cdots=t_n=0\}\subset\YY\). This curve is a chain
  comprising \(n\) copies of \(\P^1\) which generate
  \(H_2(\YY;\Z)\) as well as two copies of \(\A^1\) at either
  end of the chain. The bundle \(\MM_i\) restricts to the bundle
  \(\OO(1)\) on the \(i\)th \(\P^1\) and to the trivial bundle
  on the other \(\P^1\)s, which means that \(\LL\) restricts to
  \(\OO(1)\) on all the \(\P^1\)s. Since the compact
  irreducible components of fibres of \(\pi\colon\YY\to\YY_0\)
  are chains of \(\P^1\)s homologous to the positive linear
  combinations of \(\P^1\)s in this chain, this implies that
  \(\LL\) is relatively ample.
 
  Since the bundles \(\MM_i\) are toric line bundles generated
  by global sections, we have {\cite[Corollary on
    p.74]{FultonToric}}
  \[\OP{Ext}^j(\OO_\YY,\MM_i)=0\mbox{ for all }j>0.\] If we can
  show that \(\OP{Ext}^1(\MM_i,\OO_\YY)=0\) then we can use
  {\cite[Lemma 3.2.3]{VdB}} to deduce that
  \(\OP{Ext}^*(\OO_\YY\oplus\MM,\OO_\YY\oplus\MM)\) is supported
  in degree zero and argue as in {\cite[Proposition 3.2.5]{VdB}}
  to deduce that \(\OO_\YY\oplus\MM\) generates.

  Tensoring with $\MM_{i}^{-1}$ we see that
  $\OP{Ext}^1(\OO_\YY, \MM_i) \cong H^1(\MM_{i}^{-1})$. By
  projecting to $(t_0,\ldots, t_n)$, we can view $\mathcal{Y}$
  as a family over $\mathbb{A}^{n+1}$ which is the versal family
  of deformations of the nodal curve of the form
  $\mathbb{A}^1 \cup_{pt} \mathbb{P}^1 \cup_{pt} \mathbb{P}^1
  \cup_{pt} \ldots \mathbb{P}^1 \cup_{pt} \mathbb{A}^1$ with
  $n+1$ nodes. Any other fiber $C_t$ of this family is given by
  a nodal curve obtained from $C_0$ by smoothing the nodes
  corresponding the non-zero component of $t=(t_0,\ldots,
  t_n)$. The restriction of $\MM_i^{-1}$ to these curves gives a
  line bundle on $C_t$ whose restriction to the rational
  components of $C_t$ are either all trivial or in at most one
  component it restricts to $\mathcal{O}(-1)$. In any case,
  $H^1(\MM_{i}^{-1}|_{C_t})=0$ for any $t$, which then implies
  $H^1(\MM_{i}^{-1})=0$ as claimed.
\end{proof}

\paragraph{Corollary.} \label{cor:universal_calculation} {\em
  The derived category of \(\mathcal{Y}\) is quasi-equivalent to
  the derived category of modules over \(\AA(T^*S^1,D)\).}

\begin{proof}
  Since \(\OO_\YY \oplus \MM \) is a tilting object, the derived
  category of \(\YY\) is quasi-equivalent to the derived category
  of modules of \(\OP{End}_{\YY}(\OO_\YY \oplus \MM )\). This
  can be computed directly via toric geometry. Indeed, we have
  $\mathrm{Hom}_{\YY}(\MM_{i}, \MM_{j}) \cong H^0(\MM_{j}
  \otimes \MM_i^{-1})$ which, as in
  \ref{lma:sections_polynomials}, can be identified with the set
  of polynomials $p \in k[x_{i}, y_j]$ in the Cox ring
  such that $p(\lambda \cdot x) = \chi_{-i,j}(\lambda) p(x)$ for
  all $\lambda \in T$, where
  $\chi_{-i,j} (\lambda_1,\ldots, \lambda_n) = \lambda_i^{-1}
  \lambda_{j}$. Assuming $i \geq j$ without loss of generality, such
  polynomials are generated freely over $R$ by
  \[ x_{i} x_{i+1}\cdots x_{n} x_{0} \cdots x_{j-1} u^r,\qquad
    y_{j} y_{j+1}\cdots y_{i-1} v^s \text{\ for\ } r,s \in
    \mathbb{Z}_{\geq 0}. \] Note that
  $\mathrm{End}_{\YY}(\MM_i) \cong \mathcal{O}_\YY$
  itself is freely generated over $R$ by
  $\{1,u^r, v^s : r,s \in \mathbb{Z}_{\geq0}\}$.
  
  We now write down an isomorphism
  \(F\colon \OP{End}_\YY(\VV)\to \AA(T^*S^1,D)\). Both sides are
  bimodules over \(\bigoplus_{i=0}^n Re_i\) where \(e_i\) is an
  idempotent, acting as the identity in
  \(\Hom_\YY(\MM_i,\MM_i)= e_i\End_\YY(\VV)e_i\) or as the
  constant path at vertex \(i\) in \(\AA(T^*S^1,D)\). We define
  \[F_{ij}\colon \Hom_\YY(\MM_i,\MM_j)\to e_j\AA(T^*S^1,D)e_i\]
  using the basis above, setting
  \begin{align*}
    F_{ij}(x_{i} x_{i+1}\cdots x_{n} x_{0} \cdots x_{j-1} u^r) &= a_{j-1} \cdots a_{n} a_0 \cdots a_{i} \cdot (a_i \cdots a_n a_0 \cdots a_{i-1})^r \\
    F_{ij}(y_{j} y_{j+1}\cdots y_{i-1} v^s) &=b_j b_{j+1}
                                         \cdots b_{i-1} \cdot (b_i \cdots b_n b_0
                                         \cdots b_{i-1})^s.
  \end{align*}
  The elements on the right-hand side form a free \(R\)-module basis
  for \(e_j\AA(T^*S^1,D)e_i\), so this map is bijective. It is a
  homomorphism because it coincides with the algebra map defined by
  \(F(x_i)=a_i\), \(F(y_i)=b_i\); to see that this algebra map is
  well-defined, observe that the quiver relations
  \(a_ib_i=t_ie_{i+1}\) and \(b_ia_i=t_ie_i\) follow from
  \(x_iy_i=t_i\).
\end{proof}

One can also perform this calculation entirely within the
category of Cohen-Macaulay modules over \(\OO_{\YY_0}\); for details, see
the recent work of Zhang \cite{HaoZhang}.

\paragraph{Corollary (Base-change).} \label{cor:base_change} {\em Let \(S\) be
  a finitely generated \(R\)-algebra. Let
  \(\YY_{S,0}=\Spec(\OO_{\YY_0}\otimes_R S)\) and consider the
  diagram

\begin{center}
  \begin{tikzpicture}
    \node (A) at (0,1) {\(\YY_S\)};
    \node (B) at (2,1) {\(\YY\)};
    \node (C) at (0,0) {\(\YY _{S,0}\)};
    \node (D) at (2,0) {\(\YY_0\)};
    \node (F) at (2,-1) {\(\Spec(R)\)};
    \node (E) at (0,-1) {\(\Spec(S)\)};
    \draw[->] (A) -- (B) node [midway,above] {\(j\)};
    \draw[->] (A) -- (C);
    \draw[->] (C) -- (D);
    \draw[->] (B) -- (D);
    \draw[->] (E) -- (F) node [midway,above] {\(i\)};
    \draw[->] (D) -- (F);
    \draw[->] (C) -- (E);
    \draw[->] (A) to[out=-135,in=135] (E);
    \draw[->] (B) to[out=-45,in=45] (F);
    \node at (-1,0) {\(g\)};
    \node at (3,0) {\(f\)};
  \end{tikzpicture}
\end{center}

where \(\YY_S\) is the fibre product. The pullback \(j^*\VV\) is
a tilting bundle on \(\YY_S\) with
\(\End_{\YY_S}(j^*\VV)\cong\AA(T^*S^1,D)\otimes_R S\). In
particular, by \ref{thm:main}(B), the derived category of
perfect modules on \(\YY_S\) inherits an action of
\(\Gamma(T^*S^1,D)\).}

\begin{proof}
  The map \(\YY\to\Spec(R)\) is a conic fibration over \(\A^{n+1}\)
  with equidimensional fibres and smooth (in particular,
  Cohen-Macaulay) total space, hence flat. The endomorphism
  bundle \(\End_\YY(\VV)\) is a locally free \(\OO_\YY\)-module,
  so \(\VV\) is flat over \(\OP{Spec}(R)\) by {\cite[Lemma
    2.2]{BH}}. By {\cite[Lemma 2.9]{BH}}, this implies that
  \(j^*\VV\) is a tilting bundle with
  \(\End_{\YY_S}(j^*\VV)\cong g_*\End_{\YY_S}(j^*\VV)\cong
  i^*f_*\End_{\YY}(\VV)\cong
  i^*\AA(T^*S^1,D)\cong\AA(T^*S^1,D)\otimes_R S\). This
  base-change formula is used in the proof of {\cite[Lemma
    2.9]{BH}} but can also be found in {\cite[Lemma
    2.10]{Karmazyn}} where the pullbacks are left-derived; in
  our case all the modules are either free or locally free, so
  derived pullback equals pullback.
\end{proof}

\section{A 1-d picture of a 3-d sphere}

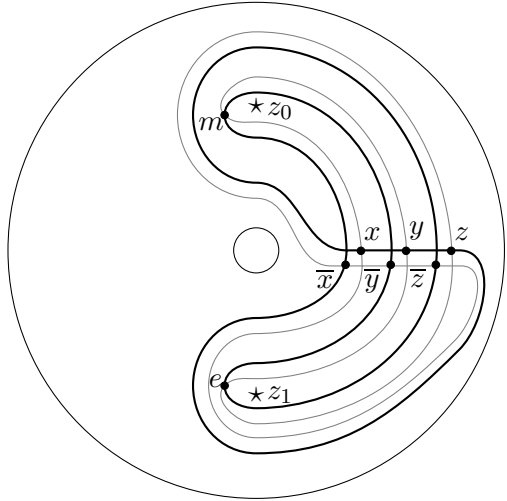
\begin{figure}[!htb]
  \begin{center}
    \begin{tikzpicture}[scale=0.6]
      \pgfmathsetmacro{\var}{0.4};
      \pgfmathsetmacro{\varh}{0.4};
      \draw (0,0) circle [radius = 5.5cm];
      \draw (0,0) circle [radius = 0.5cm];
      \node at (0,-3.2) {\(\star\)};
      \node at (0,3.2) {\(\star\)};
      \node at (0,3.1) [right] {\small \(z_0\)};
      \node at (0,-3.2) [right] {\small \(z_1\)};
      \node[minimum width=\var cm, minimum height=\varh cm] (a) at (2,0) {};
      \node[minimum width=\var cm,  minimum height=\varh cm] (b) at (3,0) {};
      \node[minimum width=\var cm, minimum height=\varh cm] (c) at (4,0) {};
      \node[minimum width=\var cm, minimum height=\varh cm] (d1) at (0,1.5) {};
      \node[minimum width=\var cm, minimum height=\varh cm] (e1) at (0,2.5) {};
      \node[minimum width=\var cm, minimum height=\varh cm] (f1) at (0,3.5) {};
      \node[minimum width=\var cm, minimum height=\varh cm] (g1) at (0,4.5) {};
      \node[minimum width=\var cm, minimum height=\varh cm] (h1) at (-0.7,3) {};
      \node[minimum width=\var cm, minimum height=\varh cm] (i1) at (-1.4,3) {}; 
      \node[minimum width=\var cm, minimum height=\varh cm] (d2) at (0,-1.5) {};
      \node[minimum width=\var cm, minimum height=\varh cm] (e2) at (0,-2.5) {};
      \node[minimum width=\var cm, minimum height=\varh cm] (f2) at (0,-3.5) {};
      \node[minimum width=\var cm, minimum height=\varh cm] (g2) at (0,-4.5) {};
      \node[minimum width=\var cm, minimum height=\varh cm] (h2) at (-0.7,-3) {};
      \node[minimum width=\var cm, minimum height=\varh cm] (i2) at (-1.4,-3) {}; 
      \draw[thick] (a.center) to[out=90,in=0] (e1.center) to[out=180,in=-90] (h1.center) to[out=90,in=180] (f1.center) to[out=0,in=90] (b.center);
      \draw[thick] (b.center) to[out=-90,in=0] (e2.center) to[out=180,in=90] (h2.center) to[out=-90,in=180] (f2.center) to[out=0,in=-90] (c.center);
      \draw[thick] (c.center) to[out=90,in=0] (g1.center) to[out=180,in=90] (i1.center) to[out=-90,in=180] (d1.center) to[out=0,in=180] (a.center);
      \draw[thick] (a.center) -- (c.center) -- (4.5,0)  to[out=0,in=45] (4.5,-2.2) to[out=-135,in=0] (g2.center) to[out=180,in=-90] (i2.center) to[out=90,in=180] (d2.center) to[out=0,in=-90] (a.center);
      \draw[gray] (a.south east) to[out=90,in=0] (e1.north) to[out=180,in=-45] (h1.center) to[out=135,in=180] (f1.north) to[out=0,in=90] (b.south east);
      \draw[gray] (b.south east) to[out=-90,in=0] (e2.south) to[out=180,in=45] (h2.center) to[out=-135,in=180] (f2.south) to[out=0,in=-90] (c.south east);
      \draw[gray] (c.south east) to[out=90,in=0] (g1.north) to[out=180,in=90] (i1.west) to[out=-90,in=180] (d1.south) to[out=0,in=180] (a.south west) -- (a.south east);
      \draw[gray] (a.south east) -- (c.south east) --++ (0.2,0)  to[out=0,in=45] (4.4,-2) to[out=-135,in=0] (g2.north) to[out=180,in=-90] (i2.east) to[out=90,in=180] (d2.south) to[out=0,in=-90] (a.south east);

      \node[shape=circle, fill=black, draw=black, inner sep=1pt] (node1) at (h1.center) {};
      \node[shape=circle, fill=black, draw=black, inner sep=1pt] (node2) at (h2.center) {};
      \node[shape=circle, fill=black, draw=black, inner sep=1pt] (node3) at (1.98,-0.32) {};
      \node[shape=circle, fill=black, draw=black, inner sep=1pt] (node4) at (2.32,-0.01) {};
      \node[shape=circle, fill=black, draw=black, inner sep=1pt] (node5) at (2.98,-0.32) {};
      \node[shape=circle, fill=black, draw=black, inner sep=1pt] (node6) at (3.32,-0.01) {};
      \node[shape=circle, fill=black, draw=black, inner sep=1pt] (node7) at (3.98,-0.32) {};
      \node[shape=circle, fill=black, draw=black, inner sep=1pt] (node8) at (4.32,-0.01) {};
     
      \node at (-1,2.8) {\small \(m\)}; 
      \node at (-0.9,-2.9) {\small \(e\)};
      \node at (1.95,-0.18) [below left] {\small \(\overline{x}\)};
      \node at (2.15,-0.01) [above right] {\small \(x\)};
      \node at (2.98,-0.16) [below left] {\small \(\overline{y}\)};
      \node at (3.15,-0.01) [above right] {\small \(y\)};
      \node at (3.98,-0.16) [below left] {\small \(\overline{z}\)};
      \node at (4.15,-0.01) [above right] {\small \(z\)};
    \end{tikzpicture}
    \caption{A 3-spherical object in \(\WW(T^*S^1,D)\) where \(|D|=2\). The gray curve is a small pushoff, used to compute the Floer complex.}\label{fig:3-sphere}
  \end{center}
\end{figure}

We conclude by discussing an example which displays how one can
draw 1-dimensional pictures corresponding to sheaves on the
higher dimensional mirrors. Let \(n=1\); in this case \(\YY\) is
the usual small-resolved conifold which is the total space of the
vector bundle $\OO(-1) \oplus \OO(-1)$ over $\P^1$. The
pushforward of the structure sheaf of $\P^1$ is well-known to be
a 3-spherical object $S$ in $D^b\Coh(\YY)$. It can be resolved
by line bundles as follows:
\[ \OO(2) \xrightarrow{(y_0,-x_1)} \OO(1)^{\oplus 2}
  \xrightarrow{(x_1,y_0)} \OO\] and $\OO(2)$ in turn is
equivalent to $\OO \xrightarrow{(x_0,y_1)} \OO(1)^{\oplus 2}$,
where $\OO(i)$ denote the line bundles on $\YY$ with degree $i$
on $\P^1$. We can, therefore, express the mirror to the
3-spherical object $S$, in terms of the generators of
$\WW(T^*S^1, D)$ and then work out, using the surgery exact
triangle on the $A$-side, which immersed Lagrangian it
corresponds to. In Figure \ref{fig:3-sphere}, the thick curve is
this immersed Lagrangian. Note that this immersed Lagrangian is
unobstructed: it does bound four ``teardrops'' (monogons) which
would contribute to the curved \(A_\infty\)-operation \(\mu_0\),
but these appear in cancelling pairs passing through the same
marked point (and hence weighted by the same variable).

The gray curve is a small pushoff. The Floer complex between
these two curves has eight generators, living in the following
degrees:

\begin{center}
\begin{tabular}{c|cccccc}
  degree & \(-2\) & \(-1\) & \(0\) & \(1\) & \(2\) & \(3\)\\
  generators& \(y\) & \(x,z\) & \(e\) & \(m\) & \(\overline{x},\overline{z}\) & \(\overline{y}\)
\end{tabular}
\end{center}

The Floer differential can be computed as follows:
\begin{align*}
  \partial y&=t_1z - t_0x, &
  \partial x &= t_1e, & \partial z &=t_0e\\
  \partial e&=0,&
  \partial m&=t_1\overline{x}-t_0\overline{z}&&\\
  \partial \overline{x}&=t_1\overline{y},& \partial \overline{z}&=t_0\overline{y},& \partial \overline{y}&=0.
\end{align*}
which yields cohomology of \(k[t_0,t_1]/(t_0,t_1)=k\) in degrees
\(0\) and \(3\).

It is also possible to verify directly that this immersed Lagrangian corresponds to a simple module of $\mathcal{A}(T^*S^1,D)$ dual to $L_0$.

\section{Derived contraction algebra}
\label{sct:dca}

\pg Let \(\YY_0\) be a 3-fold compound Du Val singularity
admitting a small resolution \(\YY\). The {\em derived
  contraction algebra} \(\Gamma\) of \(\YY\) is an enhancement
of the contraction algebra \(\Lambda\) of Donovan and Wemyss
\cite{DW16} in the sense that \(\Lambda=H^0(\Gamma)\). The
derived contraction algebra can be understood as the Drinfeld
localisation of the endomorphism algebra \(\End(\VV)\) of the
tilting bundle on \(\YY\) with respect to the idempotent
\(e=\OP{id}_{\OO_\YY}\) corresponding to the structure sheaf
\(\OO_{\YY}\). Recall that the Drinfeld localisation is given by
\[\End(\VV)_e = \End(\VV) \langle \epsilon\rangle / (\epsilon e = e\epsilon=\epsilon,\,d\epsilon =e),\]
that is we freely introduce an element \(\epsilon\) to
\(\End(\OO_\YY)\) of degree \(-1\) with \(d\epsilon=e\). This
kills the corresponding object in
\(D^b(\End(\VV))\simeq D^b(\YY)\) after localisation:
\[\mathrm{perf}(\End(\VV)_e)\simeq D^b(\YY)/\langle \OO_\YY\rangle.\]

\pg Let us consider the case of a compound $A_N$
singularity. Recall that in this case we have a 3-fold
singularity given by $uv = f_0(x,y) f_1(x,y) \cdots
f_n(x,y)$. The relative Fukaya category is derived equivalent to
the algebra $\mathcal{A}(T^*S^1, D) \otimes_R S$ where
$S\coloneqq k[x,y]$ is viewed as an $R$-algebra by the
homomorphism $t_i \to f_i(x,y)$. By Corollary
\ref{cor:base_change}, \(\mathcal{A}(T^*S^1, D) \otimes_R S\) is
isomorphic to the algebra \(\mathrm{End}_{Y_S} (j^* \VV)\) of
endomorphisms of the tilting bundle
$j^*V = \mathcal{O}_{Y_S} \oplus j^* \mathcal{M}$. Hence the
derived contraction algebra is given by
\[ \left(\AA (T^*S^1, D) \otimes_R S\right)_{e_0}, \quad
  e_0=\OP{id}_{L_0}.\] That is, the localisation of
\(D^b(\YY_S)\) away from $\mathcal{O}_{Y_S}$ corresponds to
localisation away from the Lagrangian $L_0$ in the relative
Fukaya category \(\WW(T^*S^1,D)\otimes_R S\). In the remainder
of this section, we will give an alternative, more geometric,
description of the derived contraction algebra in terms of the
relative Fukaya category of a disc.

\begin{figure}[!htb]
  \begin{center}
    \begin{tikzpicture}
      \draw[thick,->-] (10:3) arc[radius = 3,start angle=10,end angle=60];
      \draw[thick,->-] (60:3) arc[radius = 3,start angle=60,end angle=120] node [midway,above left] {\(a_{1}\)};
      \draw[thick] (120:3) arc[radius = 3,start angle=120,end angle=240];
      \draw[thick,->-] (240:3) arc[radius = 3,start angle=240,end angle=300] node [midway,below right] {\(a_{n-1}\)};
      \draw[thick,->-] (300:3) arc[radius = 3,start angle=300,end angle=350];
      \draw[thick,->-] (60:1) arc[radius = 1,start angle=60, end angle=10];
      \draw[thick,->-] (-10:1) arc[radius = 1,start angle=-10, end angle=-60];
      \draw[thick,->-] (-60:1) arc[radius = 1,start angle=-60, end angle=-120];
      \draw[thick] (-120:1) arc[radius = 1,start angle=-120, end angle=-240];
      \draw[thick,->-] (120:1) arc[radius = 1,start angle=120, end angle=60];
      \node at (-90:0.65) {\footnotesize\(b_{n-1}\)};
      \node at (90:0.6) {\footnotesize\(b_{1}\)};
      \node at (30:2cm) {\(\star\)};
      \node at (90:2cm) {\(\star\)};
      \node at (270:2cm) {\(\star\)};
      \node at (330:2cm) {\(\star\)};
      \node at (60:2.4cm) [left] {\(L_1\)};
      \node at (120:2cm) [left] {\(L_2\)};
      \node at (240:2cm) [left] {\(L_{n-1}\)};
      \node at (300:2.4cm) [right] {\(L_{n}\)};
      \node at (30:2.3cm) {\(z_0\)};
      \node at (90:2.3cm) {\(z_1\)};
      \node at (270:2.3cm) {\(z_{n-1}\)};
      \node at (330:2.3cm) {\(z_n\)};
      \node at (180:2cm) {\(\vdots\)};
      \draw[thick,->-] (10:1cm) -- (10:3cm);
      \draw[thick,->-] (-10:3cm) -- (-10:1cm);
      \draw (60:1cm) -- (60:3cm);
      \draw (120:1cm) -- (120:3cm);
      \draw (240:1cm) -- (240:3cm);
      \draw (300:1cm) -- (300:3cm);
      \node at (2,0) {excise \(L_0\)};
      \node at (3,-2) {\(T^*S^1\)};
      \node at (4,0) {\(=\)};
      \begin{scope}[shift={(8,0)}]
        \node at (3,-2) {\(\Delta\)};
        \draw (0,0) circle [radius = 3cm];
        \draw (130:3cm) -- (230:3cm);
        \draw (100:3cm) -- (260:3cm);
        \node at (0,0) {\(\cdots\)};
        \draw (80:3cm) -- (-80:3cm);
        \draw (50:3cm) -- (-50:3cm);
        \node at (-2.5,0) {\(\star\)};
        \node at (-1.2,0) {\(\star\)};
        \node at (1.2,0) {\(\star\)};
        \node at (2.5,0) {\(\star\)};
        \node at (-2.5,-0.3) {\(z_0\)};
        \node at (-1.2,-0.3) {\(z_1\)};
        \node at (1.2,-0.3) {\(z_{n-1}\)};
        \node at (2.5,-0.3) {\(z_n\)};
        \node at (-1.65,0.9) {\(L_1\)};
        \node at (-0.25,0.9) {\(L_2\)};
        \node at (1.05,0.9) {\(L_{n-1}\)};
        \node at (2.25,0.9) {\(L_n\)};
        \draw[->-,thick] (130:3) arc [radius=3,start angle=130,end angle=230] node [midway,left] {\(\alpha\)};
        \draw[->-,thick] (100:3) arc [radius=3,start angle=100,end angle=130] node [midway,above left] {\(b_1\)};
        \draw[->-,thick] (50:3) arc [radius=3,start angle=50,end angle=80] node [midway,above right] {\(b_{n-1}\)};
        \draw[->-,thick] (-80:3) arc [radius=3,start angle=-80,end angle=-50] node [midway,below right] {\(a_{n-1}\)};
        \draw[->-,thick] (-130:3) arc [radius=3,start angle=-130,end angle=-100] node [midway,below left] {\(a_1\)};
        \draw[->-,thick] (-30:3) arc [radius=3,start angle=-30,end angle=30] node [midway,right] {\(\beta\)};
      \end{scope}
    \end{tikzpicture}
  \end{center}
  \caption{Relative Fukaya category of the disc as a localisation.}
  \label{fig:rel_cat_disc}
\end{figure}
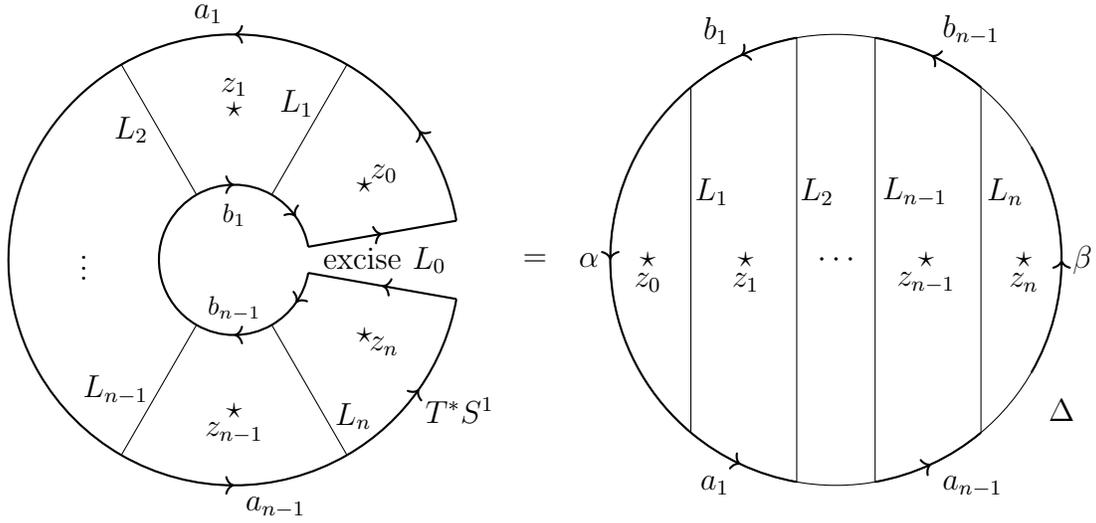

\paragraph{Theorem.} \label{thm:derived_contraction_algebra}
{\em Let \(\Delta\) be the disc obtained by excising \(L_0\)
  from \(T^*S^1\) (Figure \ref{fig:rel_cat_disc}). The derived
  contraction algebra of a 3-fold compound \(A_N\) singularity
  is quasi-equivalent to the endomorphism algebra of
  \(\bigoplus_{i=1}^n L_i\) in the relative Fukaya category
  \(\WW(\Delta,D)\otimes_R S\).}

\begin{proof}
  We can think of the annulus \(T^*S^1\) as the result of
  attaching a Weinstein 1-handle to the disc, with \(L_0\) as
  the cocore of the handle. By Ganatra, Pardon and Shende
  {\cite[Proposition 11.2]{GPS}}, this means that the
  localisation
  \[\left(\WW(T^*S^1,D)\otimes_R S\right)/\langle L_0\rangle\] is
  quasi-equivalent to the relative Fukaya category of the disc
  \(\Delta\) we get by excising \(L_0\) from \(T^*S^1\). This
  proves the theorem.
\end{proof}

\paragraph{A model for the derived contraction algebra.}
We now give a model for the \(A_\infty\)-algebra
\(\End_{\WW(\Delta,D)}\left(\bigoplus_{i=1}^nL_i\right)\). This
can be calculated directly. It is given by taking the \(R\)-linear path algebra
of the following quiver

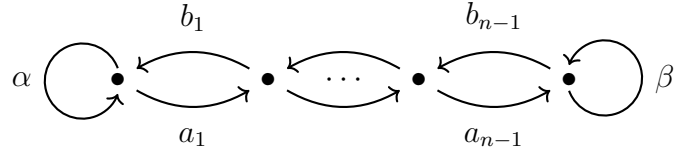
\begin{center}
  \begin{tikzpicture}
    \node (a) at (0,0) {\(\bullet\)};
    \node (b) at (2,0) {\(\bullet\)};
    \node (c) at (3,0) {\(\cdots\)};
    \node (d) at (4,0) {\(\bullet\)};
    \node (e) at (6,0) {\(\bullet\)};
    \draw[->,thick] (b) to[out=150,in=30] (a);
    \draw[->,thick] (d) to[out=150,in=30] (b);
    \draw[->,thick] (e) to[out=150,in=30] (d);
    \draw[->,thick] (d) to[out=-30,in=210] (e);
    \draw[->,thick] (b) to[out=-30,in=210] (d);
    \draw[->,thick] (a) to[out=-30,in=210] (b);
    \draw[->,thick] (0,0.15) arc [start angle=20, end angle=340,radius =0.5];
    \draw[->,thick] (6,-0.15) arc [start angle=-160, end angle=160,radius =0.5];
    \node at (-1,0) [left] {\(\alpha\)};
    \node at (7,0) [right] {\(\beta\)};
    \node at (1,0.5) [above] {\(b_1\)};
    \node at (1,-0.5) [below] {\(a_1\)};
    \node at (5,0.5) [above] {\(b_{n-1}\)};
    \node at (5,-0.5) [below] {\(a_{n-1}\)};
      \end{tikzpicture}
\end{center}

imposing the relations (coming from the quadrilaterals with
boundary \(b_i\cup L_i\cup a_i\cup L_{i+1}\) in \(\Delta\)):
\begin{gather*}
  b_i a_i = t_i e_i,\,\, a_ib_i = t_i e_{i+1},\mbox{ for
                                               }i=1,\ldots, n-1,\\
  \alpha^2=0,\qquad\beta^2=0,
\end{gather*}
and defining the differential (coming from the bigons with
boundary \(\alpha\cup L_1\) and \(L_n\cup \beta\)) by
\[ da_i=db_i=0\mbox{ for }i=1,\ldots,n-1,\qquad d \alpha = t_0
  e_1,\quad d \beta = t_n e_n,\] extending to longer paths by
the graded Leibniz rule. Note that \(a_i,b_i\),
\(i=1,\ldots,n-1\), are in degree zero whilst \(\alpha\) and
\(\beta\) are in degree \(-1\).

To see that there are no higher products, we appeal to a Maslov
index calculation of Ozsv\'{a}th and Szabo {\cite[Proposition
  6.2]{OS}} who studied these relative categories in the context
of Heegaard--Floer theory (where it is called the {\em pong
  algebra}). A rigid \((k+1)\)-gon contributing to a
\(\mu_k\)-operation has Maslov index \(2-k\); Ozsv\'{a}th and
Szabo show that the Maslov index of a holomorphic disc $u$ with
boundaries on $L_1,\ldots, L_n$ is given by
$\mathrm{mult}(u,z_1) + \mathrm{mult}(u,z_n)$, which is
non-negative since \(u\) is holomorphic. It follows that
\(k\leq 2\). A similar argument appears in {\cite[Proposition
  3.6]{Auroux}}.

\paragraph{Remark.} The relative wrapped Fukaya category $\mathcal{W}(\Delta,D)$ is acted on by its center given by its Hochschild cohomology which can be identified with the symplectic cohomology $SH(\Delta,D)$. There is a closed orbit $\eta$ that corresponds to the boundary of $\Delta$ which has degree $-2$. Thus, $\mathcal{W}(\Delta, D)$ can be seen as a category over $k[\eta]$. This recovers the familiar structure of the derived contraction algebra studied in detail in \cite[Section 6]{HuaKeller}.

\paragraph{Example.}
We can compute the case where $n=1$ and $f_0= x , f_1= y$. This
corresponds to the conifold singularity. We get that
$\Gamma= k[x,y] \langle \alpha, \beta \rangle$ with
\(\alpha^2=\beta^2=0\), \(d\alpha = x\) and \(d \beta = y\). It
is easy to determine that $H^*(\Gamma) = k [\eta]$ with
\(\eta = \alpha \beta + \beta \alpha\) of degree \(-2\). This
coincides with Booth's calculation {\cite[Section 4.2]{Booth}}.

\paragraph{Example.} Consider the Pagoda flop \(f_0 = y+x^n\),
\(f_1= y-x^n\). Our model for the derived contraction algebra
gives
\[k[x,y]\langle\alpha,\beta\rangle/(\alpha^2,\beta^2),\qquad
  d\alpha=y+x^n,\,d\beta =y-x^n.\] Assuming we are not in
characteristic \(2\), we can define
\[\zeta_1=(\alpha+\beta)/2,\qquad \zeta_2=(\alpha-\beta)/2\]
so that \(d\zeta_1=y\) and \(d\zeta_2=x^n\). This DG-algebra is
isomorphic to the graded commutative algebra
\[ k[x,y, \zeta_1, \zeta_2 ]/(\zeta_1^2+\zeta_2^2 ),\qquad d\zeta_1 = y,\, d\zeta_2 = x^n\]
Now, it is easy to see that the map from
\[k[x, \zeta],\qquad d\zeta=x^n\] sending $\zeta \to \zeta_2$ and $x \to x$ is a quasi-isomorphism.
This latter model for the derived contraction algebra of
the Pagoda flop is given by Booth in {\cite[Lemma
  4.3.8]{Booth}}. Note that in characteristic 2, the class
\(x^n\in H^0(\Gamma)\) is non-trivial, so the assumption on
characteristic is important here.
  
\paragraph{Example.} Consider the 3-fold $uv = xy
(x^2+y^3)$. This has six different partial resolutions
corresponding to different permutations of
\[f_1=x,f_2=x^2+y^3,f_3=y.\] We just focus on this particular
choice and compare the answer our model gives for
\(\Lambda = H^0(\Gamma)\) with that computed by August
{\cite[Example 4.5, Figure 2]{August}}. Our model gives
an algebra over $k[x,y]$ described by the following quiver:

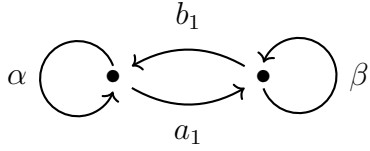
\begin{center}
  \begin{tikzpicture}
    \node (a) at (0,0) {\(\bullet\)};
    \node (b) at (2,0) {\(\bullet\)};
    \draw[->,thick] (b) to[out=150,in=30] (a);
       \draw[->,thick] (a) to[out=-30,in=210] (b);
    \draw[->,thick] (0,0.15) arc [start angle=20, end angle=340,radius =0.5];
    \draw[->,thick] (2,-0.15) arc [start angle=-160, end angle=160,radius =0.5];
    \node at (-1,0) [left] {\(\alpha\)};
    \node at (3,0) [right] {\(\beta\)};
    \node at (1,0.5) [above] {\(b_1\)};
    \node at (1,-0.5) [below] {\(a_1\)};
    \end{tikzpicture}
\end{center}

with differential \[ d\alpha = x e_1,\quad d\beta = ye_2, \]
and relations \[ a_1b_1= (x^2+y^3)e_1, \quad b_1a_1 = (x^2+y^3) e_2 , \quad \alpha^2=0, \quad \beta^2=0.\]
At the chain level, in degree zero,
we have the free \(k[x,y]\)-module spanned by
\(e_1,e_2,a_1,b_1\). We need to quotient by \[xe_1,\quad ye_2,\quad a_1b_1=(x^2+y^3)e_1,\quad
  b_1a_1=(x^2+y^3)e_2.\] The quotient algebra is therefore
generated by \(m\coloneqq ye_1\),
\(\ell\coloneqq xe_2\), \(a\coloneqq b_1\), \(c\coloneqq a_1\)
and these satisfy precisely the relations
\[\ell^2=ac,\quad m^3=ca,\quad \ell a=am=c\ell=mc=0\]
given for \(B_{con}\) in {\cite[Figure 2]{August}}. For example:
\[m^3=y^3e_1 = (y^3+x^2)e_1=a_1b_1 = ca.\]

\appendix

\section{Generation of the relative Fukaya category}
\label{sct:generation}

\renewcommand{\theparagraph}{\S A.\arabic{paragraph}\noindent}

\paragraph{Proposition} \label{prop:gen_without} {\it Let
  \(\m=(t_0,\ldots,t_n)\subset R\) and write \(k\) for the
  module \(R/\m\). The category \(\WW(T^*S^1,D)\otimes_R k\) is
  split-generated by the Lagrangian arcs \(L_0,\ldots,L_n\).}

\begin{proof}
  There is a tautological identification of
  $\mathcal{W}(T^* S^1, D )\otimes_R R/\m$ with the full
  subcategory
  \[\avoids(D)\subset \mathcal{W}(T^*S^1 \setminus D)\]
  corresponding to Lagrangian branes which do not go near the
  punctures along $D$. The manifold $T^*S^1 \setminus D$ is a
  $(n+3)$-punctured sphere with the grading structure restricted
  from the standard one on $T^*S^1$. In \cite{LP}, a mirror
  equivalence was established giving
  \[ \mathcal{W}(T^*S^1 \setminus D) \simeq D^b\Coh(C)\] where
  $C = \mathbb{A}^1 \cup_{pt} \mathbb{P}^1 \cup_{pt}
  \mathbb{P}^1 \cdots \cup_{pt} \mathbb{P}^1 \cup_{pt}
  \mathbb{A}^1$ is a nodal curve with $n+2$ irreducible toric
  components glued together at the toric fixed points. Under
  this equivalence, the full subcategory \(\avoids(D)\) gets
  identified with the full subcategory
  $\perf(C)\subset D^b\Coh(C)$, and the Lagrangians $L_i$ go to
  line bundles $\LL_i$ on $C$. In particular, one can arrange
  that \(\LL_0\) is the trivial bundle (i.e. the structure sheaf
  \(\OO_C\)).

  In the case \(n=0\), the mirror curve \(C\) is simply the {\em
    affine} curve \(\A^1\cup_{pt}\A^1=\Spec k[x,y]/(xy)\),
  and the category \(D^b\Coh(C)\) is quasi-equivalent to the
  derived category of modules over \(\OP{End}(\OO_C)\). The
  subcategory of perfect objects is then generated by
  \(\OP{End}(\OO_C)\) itself {\cite[Lemma 15.78.1]{Stacks}}.

  For higher \(n\), there is an \(n+1\)-fold covering map
  \(\pi\colon T^*S^1\setminus D\to T^*S^1\setminus\{p\}\) which
  respects the grading. The graph of \(\pi\) is a Lagrangian
  submanifold of
  \((T^*S^1\setminus D)^-\times (T^*S^1\setminus\{p\})\) (where
  \(^-\) indicates that we reverse the sign of the symplectic
  form on this factor). This induces triangulated \(A_\infty\)
  quilt functors
  \[\pi_*\colon \WW(T^*S^1\setminus D)\to \WW(T^*S^1\setminus\{p\})\quad
    \mbox{respectively}\quad\pi^*\colon \WW(T^*S^1\setminus
    \{p\})\to \WW(T^*S^1\setminus D).\] Geometrically, a
  Lagrangian brane is sent under \(\pi_*\), respectively
  \(\pi^*\), to its (possibly immersed) image, respectively
  preimage, under \(\pi\). These functors restrict to give
  functors
  \[\pi_*\colon\avoids(D)\to\avoids(p)\quad\mbox{respectively}\quad
    \pi^*\colon\avoids(p)\to\avoids(D).\] Given an object of
  \(\avoids(D)\), it follows as in {\cite[Section
    9]{SeidelGenus2}} that the object \(\pi^*\pi_*(L)\) is the
  sum \(\bigoplus_{g \in G} g(L)\) where \(G\)
  is the deck group of the covering map \(\pi\).

  Write \(L_0,\ldots,L_n\) for the arcs in \(T^*S^1\setminus D\)
  and \(\bar{L}_0\) for the arc in \(T^*S^1\setminus\{p\}\). By
  the \(n=0\) case of the proposition, if \(L\in\avoids(D)\)
  then \(\pi_*(L)\) is generated by
  \(\bar{L}_0\subset T^*S^1\setminus\{p\}\). Therefore
  \(\bigoplus_{G} g(L)\) is generated by
  \(\pi^*\bar{L}_0=\bigoplus_{i=0}^nL_i\), and since \(L\)
  is a summand of \(\bigoplus_{g \in G}g(L)\), we
  see that \(L\) is split-generated by \(\bigoplus_{i=0}^nL_i\),
  as required.
\end{proof}

\paragraph{Remark.} Obviously, the Lagrangians $L_0,\ldots, L_n$
do not generate $\mathcal{W}(T^*S^1 \setminus D)$, since the
Lagrangian branes that are allowed in
$\mathcal{W}(T^*S^1 \setminus D)$ can have ends near the
punctures along $D$.

\paragraph{Proposition (Generation with
  coefficients).}\label{prop:spectral_sequence_argument} {\it
  Let \(L\) be an object of \(\WW(T^*S^1,D)\). If \(L\) generates
  \(\WW_0(T^*S^1,D):=\WW(T^*S^1,D)\otimes_R R/\m\) then it also
  generates the relative wrapped category with coefficients in
  \(\bar{R}\), that is \(\WW(T^*S^1,D)\otimes_R\bar{R}\).}

As a corollary, the category \(\WW(T^*S^1,D)\otimes_R\bar{R}\)
is split-generated by the Lagrangian arcs
\(L_0,\ldots,L_n\). The proof of this proposition will take up
the rest of the appendix.
  
\paragraph{Proof.}
Let
\[\bar{\AA}=\OP{End}_{\WW(T^*S^1,D)}(L)\otimes_R\bar{R},\qquad\AA_0=\OP{End}_{\WW_0(T^*S^1,D)}(L)
  = \OP{End}_{\WW(T^*S^1,D)}(L)\otimes_R R/\m.\]

We have Yoneda-type functors
\[\bar{\Yon}\colon\WW(T^*S^1,D)\otimes_R\bar{R}\to
  \mod(\bar{\AA})\]
and \[\Yon_0\colon\WW_0(T^*S^1,D)\to\mod(\AA_0).\] The module
\(\Yon_0(L)=\AA_0\) (respectively \(\bar{\Yon}(L)=\bar{\AA}\))
generates the subcategory \(\perf(\AA_0)\) (respectively
\(\perf(\bar{\AA})\)) of perfect objects. Since \(L\) generates
\(\WW_0(T^*S^1,D)\), the functor \(\Yon_0\) lands in
\(\perf(\AA_0)\) and corestricts to give a quasi-equivalence
  \[\Yon_0\colon \WW_0(T^*S^1,D)\to\perf(\AA_0)\]
  (i.e. the induced functor on homotopy categories is fully
  faithful and essentially surjective). We want to show that
  \begin{enumerate}
  \item[(a)] \(\bar{\Yon}\) lands in \(\perf(\bar{\AA})\);
  \item[(b)] the induced functor \(H(\bar{\Yon})\) on homotopy
    categories is (i) essentially surjective and (ii) fully
    faithful.
  \end{enumerate}

\paragraph{Proof of (a):} The subcategory
  \(\perf(\bar{\AA})\subset\mod(\bar{\AA})\) is precisely the
  triangulated subcategory of compact objects (see for example
  {\cite[Proposition 15.78.3]{Stacks}}). An object \(C\) in a
  pre-triangulated $A_\infty$ category is compact if and only if the functor it corepresents $\mathrm{hom}(C,\,\ )$ preserves coproducts, that is,  
  \[ \oplus_i \hom(C,E_i)=\hom(C,\oplus_i E_i) \] for
  arbitrary direct sums \(\oplus_iE_i\). So it suffices to
  show that if \(K\in \WW(T^*S^1,D)\otimes_R\bar{R}\) is an
  object then
  \[\oplus_i\hom_{\mod(\bar{\AA})}\left(\bar{\Yon}(K),E_i\right) =
    \hom_{\mod(\bar{\AA})}\left(\bar{\Yon}(K),\oplus_i
      E_i\right)\] for arbitrary direct sums \(\oplus_iE_i\)
  in \(\mod(\bar{\AA})\).

  The complexes
  \(\oplus_i\hom_{\mod(\bar{\AA})}\left(\bar{\Yon}(K),E_i\right)\)
  and
  \(\hom_{\mod(\bar{\AA})}\left(\bar{\Yon}(K),\oplus_i
    E_i\right)\) are complete filtered \(\bar{R}\)-modules with
  the filtration coming from the action of powers of the maximal
  ideal; the canonical map
  \begin{equation}\label{eq:canonical_map}
    \oplus_i\hom_{\mod(\bar{\AA})}\left(\bar{\Yon}(K),E_i\right) \to
    \hom_{\mod(\bar{\AA})}\left(\bar{\Yon}(K),\oplus_i E_i\right)
  \end{equation} is a morphism of
  filtered complexes. There are therefore spectral sequences
  computing both sides, and a morphism of spectral sequences
  induced by \eqref{eq:canonical_map}. By the Eilenberg-Moore
  comparison theorem, it suffices to check that this morphism is
  an isomorphism on the \(E_0\) pages. Note that Eilenberg-Moore
  requires completeness of the filtration, which is why we are
  working over \(\bar{R}\) instead of \(R\).

  The \(E_0\) pages are respectively
  \begin{align*}
    E_0^{pq} = \oplus_i\hom^{p+q}_{\mod(\AA_0)}\left(\Yon_0(K),\OP{gr}^p(E_i)\right) \mbox{
      and } E_0^{pq} = \hom_{\mod(\AA_0)}^{p+q}\left(\Yon_0(K),\oplus_i
      \OP{gr}^p(E_i)\right)
  \end{align*} where \(\OP{gr}^p\) denotes the \(p\)th graded
  piece of the associated graded module. The morphism on
  \(E_0\)-pages is induced by the canonical map
  \[\oplus_i\hom_{\mod(\AA_0)}\left(\Yon_0(K),\OP{gr}(E_i)\right)\to
    \hom_{\mod(\AA_0)}\left(\Yon_0(K),\oplus_i
      \OP{gr}(E_i)\right).\] Since \(\Yon_0(K)\) is perfect, this
  is an isomorphism, which proves (a).

\paragraph{Proof of (b.i):} We have \(\bar{\AA}=\bar{\Yon}(L)\), and
  since \(\bar{\AA}\) generates \(\perf(\bar{\AA})\), the
  essential image of \(\bar{\Yon}\) in \(\mod(\bar{\AA})\)
  contains \(\perf(\bar{\AA})\).

\paragraph{Proof of (b.ii):} Given objects
  \(K,K'\in\WW(T^*S^1,D)\otimes_R\bar{R}\), the complexes
  \[CF(K,K')\otimes_R\bar{R}\quad\mbox{ and
    }\quad\hom_{\mod(\bar{\AA})}\left(\bar{\Yon}(K),\bar{\Yon}(K')\right)\]
  are filtered by powers of the maximal ideal. These filtrations
  give us spectral sequences and the functor \(\bar{\Yon}\) gives
  a map of filtered complexes
  \(CF(K,K')\otimes_R\bar{R}\to
  \hom_{\mod(\bar{\AA})}\left(\bar{\Yon}(K),\bar{\Yon}(K')\right)\)
  and hence a morphism of spectral sequences. On the \(E_1\)
  page this is just the map
  \[H\left(\hom_{\WW_0(T^*S^1,D)}(K,K')\right)\otimes_R
    \OP{gr}(\bar{R})\to
    H\left(\hom_{\mod(\AA_0)}(\Yon_0(K),\Yon_0(K'))\right)\otimes_R
    \OP{gr}(\bar{R})\] induced from
  \(H(\Yon_0)\colon H\left(\hom_{\WW_0(T^*S^1,D)}(K,K')\right)\to
  H\left(\hom_{\mod(\AA_0)}(\Yon_0(K),\Yon_0(K'))\right)\)
  (because any polygons which pass through the marked points
  have their contributions weighted by an element of
  \(\m\)). This is an isomorphism because \(\Yon_0\) is
  cohomologically full and faithful. The Eilenberg-Moore
  comparison theorem then implies that the map \(H(\bar{\Yon})\)
  is an isomorphism, which proves that \(\bar{\Yon}\) is
  cohomologically full and faithful.\qed

\Addresses

\end{document}